%% file: paperfireshape.tex
\RequirePackage{fix-cm}
\documentclass[twocolumn]{svjour3}          
\def\makeheadbox{\relax}
\smartqed  
\usepackage{graphicx}
\usepackage{amsmath, amssymb}
\usepackage{latexsym}
\usepackage{mathtools}
\usepackage{color}
\usepackage{cite}
\usepackage{tikz}
\usepackage{url}
\usepackage{newfloat}
\DeclareFloatingEnvironment[name=Lst.]{mylisting}

\include{slidesymbols}

\begin{document}

\title{Fireshape: a shape optimization toolbox for Firedrake}

\author{Alberto Paganini         \and
        Florian Wechsung 
}

\institute{A. Paganini \at
              School of Mathematics and Actuarial Science, University of Leicester, University Road, Leicester, LE1 7RH, United Kingdom, \\
              \email{a.paganini@leicester.ac.uk}
           \and
           F. Wechsung \at
              Courant Institute of Mathematical Sciences, New York University, 251 Mercer St, New York, NY 10012,\\
              \email{wechsung@nyu.edu}
}

\date{}

\maketitle

\begin{abstract}
We introduce Fireshape, an open-source and automated shape optimization toolbox for the
finite element software Firedrake. Fireshape is based on the moving mesh
method and allows users with minimal shape optimization knowledge
to tackle with ease challenging shape optimization problems constrained to partial
differential equations (PDEs).
\end{abstract}

\section{Introduction}
\label{sec:intro}

One of the ultimate goals of structural optimization is the development of
fully automated software that allows users to tackle challenging structural
optimization problems in the automotive, naval, and aerospace industries
without requiring deep knowledge of structural optimization theory. The
scientific community is working actively in this direction, and recent years
have seen the publication of educational material that simplifies the
understanding of structural optimization algorithms and guides the
development of related optimization software. These resources are based on
different models, such as moving mesh methods \cite{AlPa06,DaFrOm18,
CoLePrFrLa19, Al07}, level-sets \cite{La18, AlJoTo02, BeWaBe19}, phase fields
\cite{GaHeHiKa15, DoPoRuSi20, BlGaFaHaSt14}, and SIMP\footnote{The acronym
SIMP stands for Solid Isotropic Material with Penalisation.} \cite{Si01,
polytop12, Be89, BeSi04}, and are implemented in various software environments such
as Matlab \cite{Si01, polytop12}, FreeFem++ \cite{AlPa06,DaFrOm18,AlJoTo02},
OpenFOAM CFD \cite{CoLePrFrLa19}, FEMLAB \cite{LiKoHu05}, and FEniCS
\cite{La18, DoFuJoSc19}, to mention just a few.

In this work, we introduce Fireshape: an automated shape optimization library
based on the moving mesh approach that requires very limited input from the
user. Shape optimization refers to the optimization of domain boundaries and
plays an important role in structural design. For instance, shape
optimization plays a crucial role in the design of airfoils
\cite{HeJiMaYiMa19, WaDoJi18, ScIlScGa13} and boat hulls
\cite{LoPaQuRo12,MaVuCu16,AbSu17}. Shape optimization is also a useful
refinement step to be employed after topology optimization \cite[Ch.
1.4]{BeSi04}. Indeed, topology optimization allows more flexibility in
geometric changes and it is a powerful tool to explore a large design space.
However, topology optimization {\color{black} may return} slightly blurred
(grey-scale) and/or staircase designs \cite{BeSi04}. By adding a final shape
optimization step, it is possible to post-process results computed with
topology optimization and devise optimal designs with sharp boundaries and
interfaces, {\color{black} when this is necessary}.

Fireshape is based on the moving mesh shape optimization approach \cite[Ch.
6]{Al07}. In this approach, geometries are parametrized with meshes that can
be arbitrarily precise and possibly curvilinear. The mesh nodes and faces are
then optimized (or ``moved'') to minimize a chosen target function. 
Fireshape has been developed on the moving mesh approach
because the latter has a very neat interpretation in terms of geometric
transformations and is inherently compatible with standard finite element
software, {\color{black} as it relies on the canonical construction of
finite elements via pullbacks} \cite{PaWeFa18}. The main drawback of the moving mesh
approach is that it does not allow topological changes in a straightforward
and consistent fashion. However, Fireshape has been developed to facilitate
shape optimization, and topology optimization is beyond its scope.

Fireshape couples the finite element library Firedrake \cite{Firedrake1,
Firedrake2, Firedrake3, Firedrake4, Firedrake5} with the Rapid Optimization
Library (ROL) \cite{ROL}. Fireshape allows decoupled discretization of
geometries and state constraints and it includes all necessary routines to
perform shape optimization (geometry updates, regularization terms, geometric
constraints, etc.). To solve a shape optimization problem in Fireshape, users
must describe the objective function and the eventual constraints using the
Unified Form Language (UFL) \cite{AlLoOlRoWe14}, a {\color{black}language
embedded in Python to describe variational forms and their finite element
discretizations} that is very similar to standard mathematical notation. Once
objective functions and constraints have been implemented with UFL, users
need only to provide a mesh that describes the initial design and, finally,
select their favorite optimization algorithm from the optimization library
ROL{\color{black}, which contains algorithms for unconstrained, bound
constrained, and (in-)equality constrained optimisation problems. In
particular, users do not need to spend time deriving shape derivative
formulas and adjoint equations by hand, nor worry about their discretization,
because Fireshape automates these tasks.}

Typically the {\color{black} computational} bottleneck of PDE constrained
optimization code lies in the solution of the state and the adjoint equation.
While Fireshape and Firedrake are both written in Python, to assemble the
state and adjoint equations the Firedrake library automatically generates
optimized kernels in the C programming language. The generated systems of
equations are then passed to the PETSc library (also written in C) and can
be solved using any of the many linear solvers and preconditioners accessible from
PETSc \cite{petsc1, petsc2, petsc3, petsc4, petsc5, petsc6}. This combination
of Python for user facing code and C for performance critical parts is well
established in scientific computing as it provides highly performant code
that is straightforward to develop and use. Finally, we mention that
Fireshape, just as Firedrake, PETSc and ROL, supports MPI parallelization and
hence can be used to solve even large scale three dimensional shape
optimization problems. {\color{black}In particular, Firedrake has been used
to solve systems with several billion degrees of freedom on supercomputers
with tens of thousands of cores (see e.g.~\cite{mitchell_high_2016,
farrell2019augmented}), and we emphasize that any PDE solver written in
Firedrake can be used for shape optimization problems with Fireshape.}

Fireshape is an open-source software licensed under the GNU Lesser General
Public License v3.0. Fireshape can be downloaded at
\begin{center}
\url{https://github.com/Fireshape/Fireshape} .
\end{center}
Its documentation contains
several tutorials and is available at
\begin{center}
\url{https://fireshape.readthedocs.io/en/latest/} .
\end{center} To illustrate Fireshape's
capabilities and ease of use, in Section \ref{sec:pipe} we provide a tutorial
and solve a three-dimensional shape optimization problem constrained by the
nonlinear Navier-Stokes equations. The shape optimization knowledge required
to understand this tutorial is minimal. {\color{black}In Section
\ref{sec:cantilever}, we provide another tutorial based on a shape
optimization problem constrained to the linear elasticity equation.} In Section
\ref{sec:shapeoptimization}, we describe in detail the rich mathematical
theory that underpins Fireshape. In Section \ref{sec:anatomy}, we describe
Fireshape's main classes and Fireshape's extended functionality. Finally, in
Section \ref{sec:conclusions}, we provide concluding remarks.

\section{Example: Shape optimization of a pipe}
\label{sec:pipe}

\begin{figure}[htb!]
\center
\begin{tikzpicture}[scale=0.6]
\draw[thick, yshift=1cm] (0,0)--(1.5,0);--(6,2)--(9,2);
\draw[thick] (0,0)--(1.5,0);
\draw[thick] (1.5,0) to[out=0,in=210] (3,0.2);
\draw[thick] (3,0.2) to[out=30,in=210] (6,0.2+3*0.5773);
\draw[thick] (6,0.2+3*0.5773) to[out=30,in=180] (7.5,0.4+3*0.5773);
\draw[thick] (7.5,0.4+3*0.5773) -- (9,0.4+3*0.5773);
\draw[thick, yshift=1cm] (1.5,0) to[out=0,in=210] (3,0.2);
\draw[thick, yshift=1cm] (3,0.2) to[out=30,in=210] (6,0.2+3*0.5773);
\draw[thick, yshift=1cm] (6,0.2+3*0.5773) to[out=30,in=180] (7.5,0.4+3*0.5773);
\draw[thick, yshift=1cm] (7.5,0.4+3*0.5773) -- (9,0.4+3*0.5773);
\draw[very thick, ->, xshift =-1cm] (-0.5,0.5) -- (0.5,0.5);
\draw[very thick, ->, xshift=10cm, yshift=0.4cm+3*0.5773cm] (-0.5,0.5) -- (0.5,0.5);
\draw[thick] (0,0.5) ellipse (0.3cm and 0.5cm);
\draw[thick, dashed, yshift = 0.9cm+3*0.5773cm, xshift = 9cm] (270:0.3cm and 0.5cm) arc (270:90:0.3cm and 0.5cm);
\draw[thick, yshift = 0.9cm+3*0.5773cm, xshift = 9cm] (-90:0.3cm and 0.5cm) arc (-90:90:0.3cm and 0.5cm);
\end{tikzpicture}
\caption{Viscous fluid flows through a pipe $\rm{\Omega}$
from the left to the right. A poor pipe design can lead to 
an excessive amount of kinetic energy being dissipated into heat.}
\label{fig:sketch}
\end{figure}
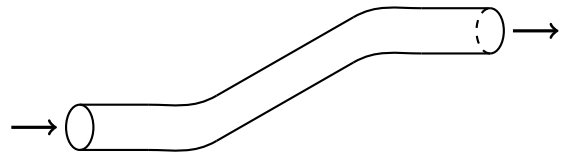

For this hands-on introduction to Fireshape, we {\color{black}first}
consider a viscous fluid flowing through a pipe $\rm{\Omega}$ (see Figure
\ref{fig:sketch}) and aim to minimize the kinetic energy dissipation into
heat by optimizing the design of $\rm{\Omega}$. {\color{black} This is a
standard test case \cite{ZhLi08,ViMa17,FeAlDaJo20, DiDiFuSiLa18,DaFaMiAl17}.}
To begin with, we need to describe this optimization problem with a
mathematical model.

We assume that the fluid velocity $\Vu$ and the fluid pressure $\rm p$
satisfy the incompressible Navier-Stokes equations  \cite[Eqn. 8.1]{ElSiWa14}, which read

\begin{subequations}
\begin{align}
\label{eq:NSmain}
    -2\nu \nabla \cdot \varepsilon\Vu + (\Vu\cdot \nabla) \Vu + \nabla \rm{p} &= 0 & \text{in } \rm{\Omega}\,,\\
\label{eq:NSincompressible}
\Div \Vu &= 0& \text{in } \rm{\Omega}\,,\\
\label{eq:NSdirBC}
\Vu &= \Vg &\text{on } \partial\rm{\Omega}\setminus \rm{\Gamma}\,,\\
\label{eq:NSrobBC}
\rm{p}\Vn - 2\nu \varepsilon \Vu \, \Vn & = 0 &\text{on }  \rm{\Gamma}\,,
\end{align}
\label{eq:NavierStokes}
\end{subequations}
where $\nu$ denotes the fluid viscosity, $\varepsilon \Vu = \varepsilon (\Vu)
= \frac{1}{2}(\nabla \Vu+\nabla\Vu^\top)$, $\nabla \Vu$ is the derivative
(Jacobian matrix) of $\Vu$ and $\nabla\Vu^\top$ is the derivative transposed,
$\rm{\Gamma}$ denotes the outlet, and the function $\Vg$ is a Poiseuille flow
velocity \cite[p. 122]{ElSiWa14} at the inlet and zero on the pipe walls.
In our numerical experiment, the inlet is a disc of radius 0.5 in the
xy-plane, and $\Vg(x,y,z) = (0,0,1 - 4(x^2 +y^2))^\top$ on the inlet.

To model the kinetic energy dissipation, we consider the diffusion term
in \eqref{eq:NSmain} and introduce the function
\begin{equation}\label{eq:kineticfunctional}
\rm{J}(\rm{\Omega}) = \int_{\rm{\Omega}} \nu\, \varepsilon \Vu : \varepsilon\Vu \DX\,,
\end{equation}
where the
colon symbol $:$ denotes the Frobenius inner product, that is,
$\varepsilon \Vu : \varepsilon \Vu= \rm{trace}((\varepsilon \Vu)^\top \varepsilon \Vu)$.

Now, we can formulate the shape optimizations problem we are considering.
This reads:
\begin{equation}\label{eq:shapeoptproblem}
\text{Find } {\rm \Omega}^*
\text{ such that } \rm{J}(\rm{\Omega}^*) = \inf \rm{J}(\rm{\Omega}) \text{ subject to
\eqref{eq:NavierStokes}.}
\end{equation}
To make this test case more interesting, we further impose a volume equality
constraint on $\rm{\Omega}$. Otherwise, {\color{black} without this volume
constraint}, the solution to this problem would be a pipe with arbitrarily
large diameter.

In the next subsections, we explain step by step how to solve this shape
optimization problem in Fireshape. To this end, we need to create: a mesh
that approximates the initial guess $\rm \Omega$ (Section \ref{sec:step1}),
an object that describes the PDE-constraint \eqref{eq:NavierStokes} (Section
\ref{sec:step2}), an object that describes the objective function
\eqref{eq:kineticfunctional} (Section \ref{sec:step3}), and, finally, a
``main file'' to set up and solve the optimization problem
\eqref{eq:shapeoptproblem} {\color{black}with the additional volume
constraint} (Section \ref{sec:finalstep}).

\subsection{Step 1: Provide an initial guess}
\label{sec:step1}
The first step is to provide a mesh that describes an initial guess of $\rm
\Omega$. For this tutorial, we create this mesh using the software Gmsh
\cite{GeRe09}. The initial guess employed is sketched in Figure
\ref{fig:sketch}. 
{\color{black} The only mesh detail necessary to understand this
tutorial is that, in Listing \ref{lst:PDEconstraint},
the number 10 corresponds to the inlet, wheres the
boundary flags 12 and 13 refer to the pipe's boundary.}
For the remaining geometric details, we refer to the code archived on
Zenodo~\cite{zenodo-firedrake-and-driver}.

\subsection{Step 2: Implement the PDE-constraint}
\label{sec:step2}
The second step is to implement a finite element solver
of the Navier-Stokes equations \eqref{eq:NavierStokes}.
To derive the weak formulation of \eqref{eq:NSmain},
we multiply
Equation \eqref{eq:NSmain} with a test (velocity) function $\Vv$ that vanishes
on $\partial\rm{\Omega}\setminus \rm{\Gamma}$,
and Equation \eqref{eq:NSincompressible} with a test (pressure) function $\rm{q}$.
Then, we integrate over $\rm \Omega$,
integrate by parts \cite[Eq. 3.18]{ElSiWa14}, and impose the boundary condition
\eqref{eq:NSrobBC}. The resulting weak formulation of equations
\eqref{eq:NavierStokes} reads:
\begin{align}\label{eq:NSweak}
\nonumber&\text{Find }\Vu\text{ such that  }\Vu = \Vg\text{ on }
\partial\rm \Omega\setminus \Gamma\text{ and}\\
\nonumber&\int_{\rm{\Omega}} 2\nu \varepsilon \Vu :\varepsilon \Vv -  \rm{p}\Div(\Vv)
+ \Vv\cdot (\nabla \Vu)\Vu + \rm{q}\Div(\Vu) \DX = 0\\
&\text{for any pair }(\Vv, \rm{q})\,.
\end{align}

\begin{mylisting}[htb!]
\includegraphics[width=\linewidth]{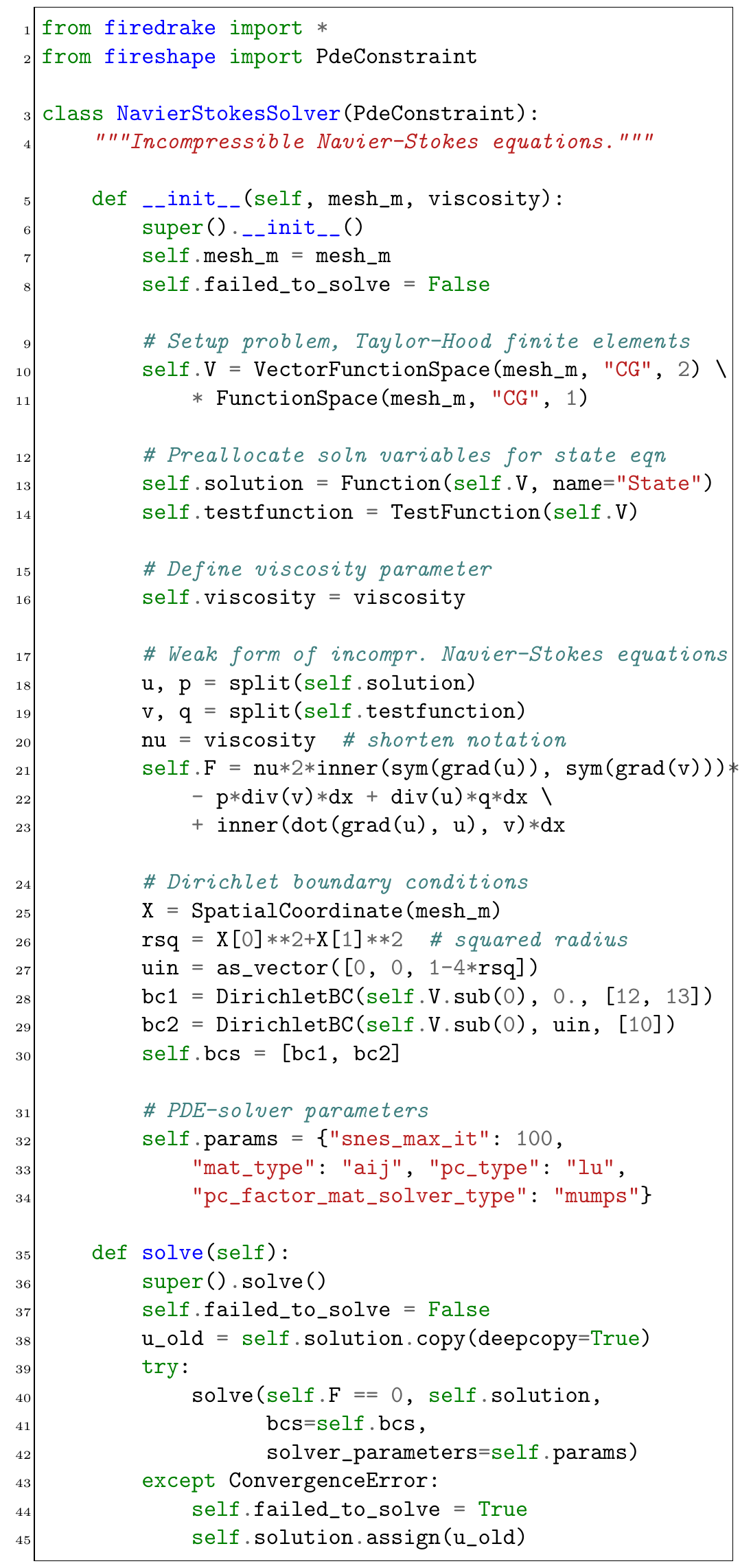}
\caption{Incompressible Navier-Stokes solver
based on Taylor-Hood finite elements.}
\label{lst:PDEconstraint}
\end{mylisting}

To implement the finite element discretization of \eqref{eq:NSweak}, we
create a class \verb!NavierStokesSolver! that inherits from the Fireshape's
class \verb!PdeConstraint!, see Listing \ref{lst:PDEconstraint}. To
discretize \eqref{eq:NSweak}, we employ P2-P1 Taylor-Hood finite elements,
that is, we discretize the trial and test velocity functions $\Vu$ and $\Vv$
with piecewise quadratic Lagrangian finite elements, and the trial and test
pressure functions $\rm{p}$ and $\rm{q}$ with piecewise affine Lagrangian
finite elements. It is well known that this is a stable discretization of the
incompressible Navier-Stokes equations \cite[pp 136-137]{ElSiWa14}.

To address the nonlinearity in \eqref{eq:NSweak}, we use PETSc's Scalable
Nonlinear Equations Solvers (SNES) \cite{petsc1, petsc2, petsc3, petsc4,
petsc5, petsc6}. In each iteration, SNES linearizes Equation
\eqref{eq:NSweak} and solves the resulting system with a direct solver. In
general, it is possible that the SNES solver fails to converge sufficiently
quickly (in our code, we allow at most 100 {\color{black}SNES} iterations).
For instance, this happens if the finite element mesh self-intersects, or if
the initial guess used to solve \eqref{eq:NSweak} is not sufficiently good.
Most often, these situations happen when the optimization algorithm takes an
optimization step that is too large. In these cases, we can address SNES'
failure to solve \eqref{eq:NSweak} by reducing the optimization step. In
practice, we deal with these situations with Python's \texttt{try: ...
except: ...} block. If the SNES solver fails with a
\texttt{ConvergenceError}, we catch this error and set the boolean flag
\texttt{self.failed\_to\_solve} to \texttt{True}. This flag is used to adjust
the output of the objective function $\rm{J}$, see Section \ref{sec:step3}.

\subsection{Step 3: Implement the objective function}
\label{sec:step3}

The third step is to implement a code to evaluate the objective function $\rm{J}$
defined in Equation \eqref{eq:kineticfunctional}.
For this, we create a class \texttt{EnergyDissipation} that inherits
from Fireshape's class \texttt{ShapeObjective},
see Listing \ref{lst:objective}. Of course, evaluating
$\rm{J}$ requires access to the fluid velocity $\Vu$.
This access is implemented by assigning the
variable \texttt{self.pde\_solver}. 
This variable gives access also to
the variable \texttt{NavierStokesSolver.failed\_to\_solve},
which can be used to control the output of $\rm{J}$ when the
Navier-Stokes' solver fails to converge. Here, we decide that the value of
$\rm{J}$ is \texttt{NaN} (``not a number'')
if the Navier-Stokes' solver fails to converge.

\begin{mylisting}[htb!]
\includegraphics[width=\linewidth]{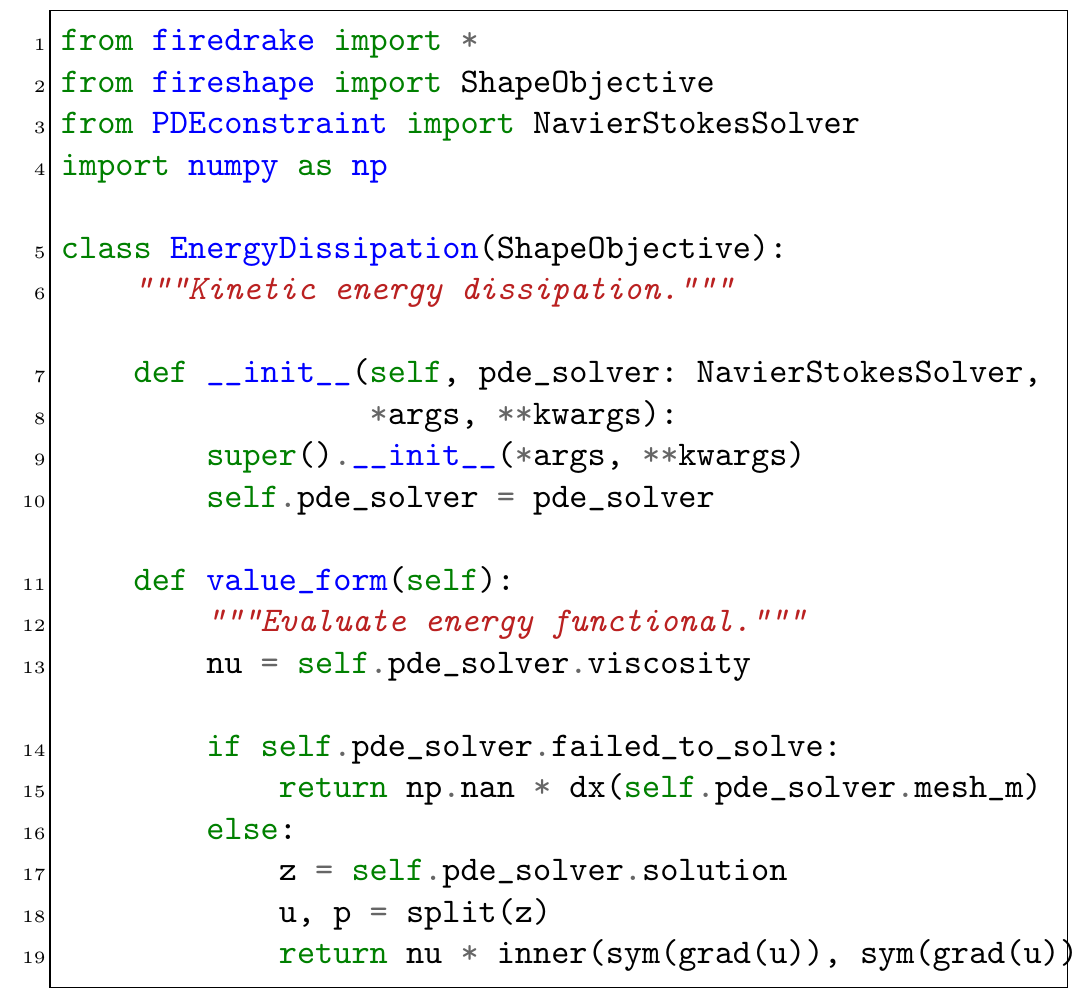}
\caption{Objective function $\rm J$ \eqref{eq:kineticfunctional},
which quantifies the kinetic energy
dissipation of the fluid. Note that $\rm J$ returns \texttt{NaN} (``not a number'')
if the solver of the PDE-constraint \eqref{eq:NSweak}
fails to converge.}
\label{lst:objective}
\end{mylisting}

\subsection{Final step: Set up and solve the problem}
\label{sec:finalstep}
At this stage, we have all the necessary ingredients to tackle
the shape optimization problem \eqref{eq:shapeoptproblem}.
The final step is to create a ``main file'' that loads the initial mesh, sets-up
the optimization problem, and solves it with an optimization algorithm.
The rest of this section contains a line-by-line description of
the ``main file'', which is listed in Listing \ref{lst:main}.
\begin{mylisting}[htb!]
\includegraphics[width=\linewidth]{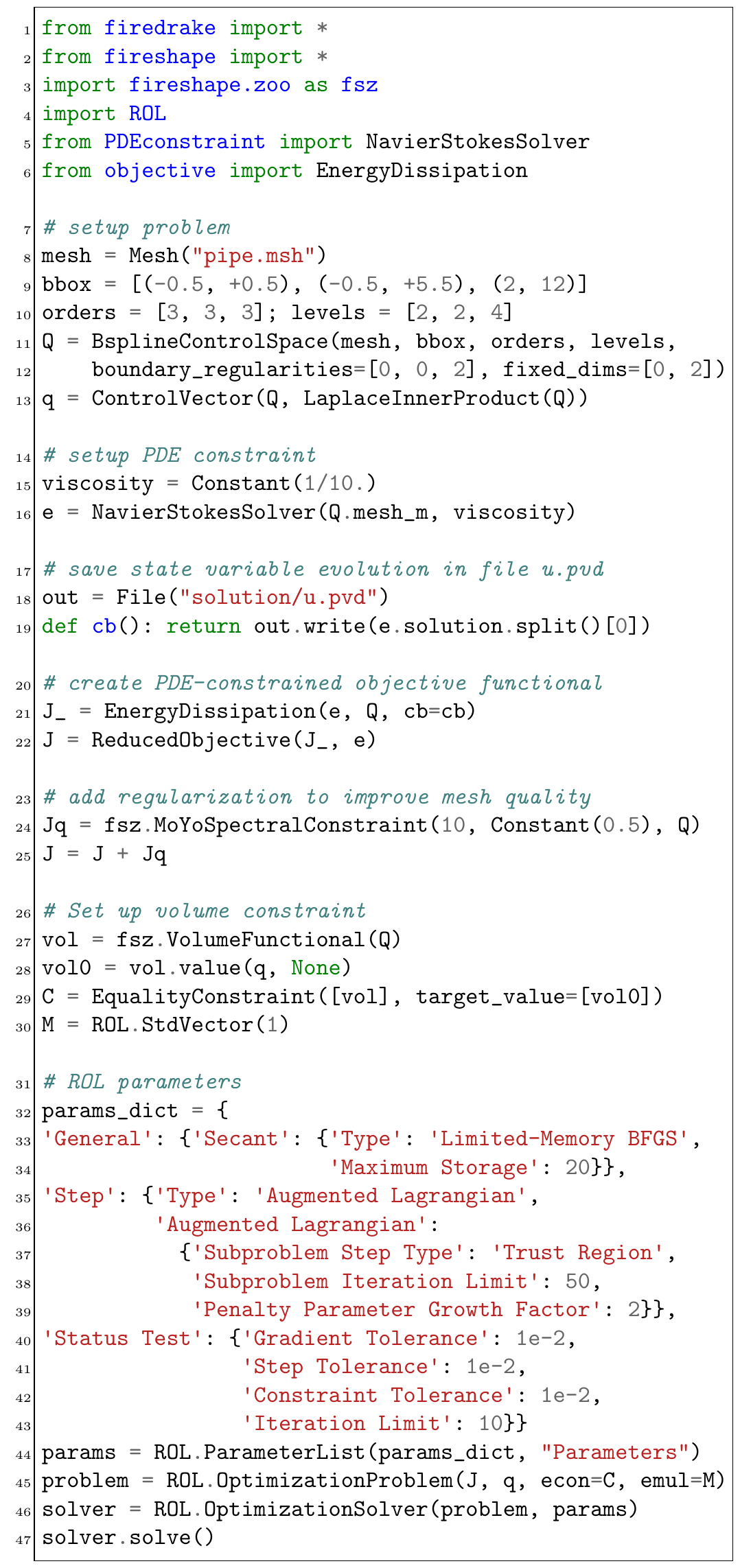}
\caption{Minimize the kinetic energy dissipation of a fluid using Fireshape.}
\label{lst:main}
\end{mylisting}
\renewcommand\labelitemi{\textbullet}
\begin{itemize}
\item \emph{In lines 1-6,} we import all necessary Python libraries and modules.
\item \emph{In lines 8-13,} we load the mesh of the initial design and we
specify {\color{black} to discretize perturbations of the domain $\rm
\Omega$ using quadratic B-splines with support in a box that does not include
the inlet and the outlet. In this example, we also impose some boundary
regularity on B-splines and limit perturbations to be only in direction of the $y$-axis.} Sections
\ref{sec:controlspace} and \ref{sec:innerproduct} \textcolor{black}{provide}
more details about the Fireshape classes \texttt{ControlSpace} and
\texttt{InnerProduct}.
\item \emph{In lines 15-16,} we initiate the Navier-Stokes finite element solver.

\item \emph{In lines 18-19,} we tell Fireshape to store the finite element solution
to Navier-Stokes' equations in the folder \texttt{solution} every time the domain
$\rm \Omega$ is updated. We can visualize the evolution of
$\Vu$ along the optimization process by opening
the file \texttt{u.pvd} with Paraview \cite{ahrens2005paraview}.

\item \emph{In lines 21-22,} we initiate the objective function $\rm J$ and
the associated reduced functional, which is used by Fireshape to define the
appropriate Lagrange functionals to shape differentiate $\rm J$. We refer to
Section \ref{sec:shapeoptshapecalc} for more details about shape
differentiation using Lagrangians. We stress that Fireshape does not require
users to derive (and implement) shape \textcolor{black}{derivative} formulas
by hand. The whole shape differentiation process is automated using pyadjoint
\cite{HaMiPaWe19, DoMiFu20}{\color{black}, see Remark \ref{rmk:shapediff}}.

\item \emph{In lines 24-25,} we add an additional regularization term to $\rm
J$ to promote mesh quality of domain updates. This regularization term
controls the pointwise maximum singular value of the geometric transformation
used to update $\rm \Omega$ \cite{wechsungthesis}. See Section
\ref{sec:shapeoptshapecalc}, Figure \ref{fig:Q}, \textcolor{black}{and
Remark \ref{rmk:detDTsmall}} for more information about the role of geometric
transformations in Fireshape.

\item \emph{In lines 27-30,} we set up an equality constraint
to ensure that the volume of the initial and the optimized domains
are equal.

\item \emph{In lines 32-43,} we select an optimization algorithm
from the optimization library ROL.
More specifically, we use an augmented Lagrangian approach
\cite[Ch. 17.3]{NoWr06} with limited-memory BFGS Hessian updates
\cite[Ch. 6.1]{NoWr06} to deal with the volume equality constraint,
and a Trust-Region algorithm \cite[Ch. 4.1]{NoWr06} to solve the intermediate
models generated by the augmented Lagrangian algorithm.
\item \emph{In lines 44-47,} finally, we gather all information and solve the problem.
\end{itemize}

{\color{black} \begin{remark} \label{rmk:shapediff} A reader may wonder why
Fireshape does not require information about shape derivative formulas and
adjoint equations to solve a PDE-constrained shape optimization problem. This
is possible because Fireshape employs UFL's and pyadjoint's automated shape
differentiation and adjoint derivation capabilities \cite{HaMiPaWe19,
DoMiFu20}. In particular, UFL combines the automated construction of finite
element pullbacks with symbolic differentiation to automate the evaluation of
directional shape derivatives along vector fields discretized by finite
elements \cite{HaMiPaWe19}. This process mimics the analytical
differentiation of shape functions and is equivalent to the
``optimize-then-discretize'' paradigm, i.e.~it yields exact (up to floating
point accuracy) derivatives. \end{remark}}

\subsection{Results}
\label{sec:results}

Running the code contained in Listing~\ref{lst:main} optimizes the pipe design when the fluid viscosity $\nu$ is 0.1, which corresponds to Reynold-number $\Re = 1./\nu=10$.
Of course, the resulting shape depends on the fluid viscosity.
A natural question is how the optimized shape depends on this parameter.
To answer this question, we can simply run Listing~\ref{lst:main} for different Reynold-numbers (by modifying line 13 with the desired value) and inspect the results.
Here we perform this comparison for $\Re \in \{1, 250, 500, 750, 1000\}$.\footnote{For this computation we extend the code in Listing~\ref{lst:PDEconstraint} and include
analytic continuation in the Reynolds-number to re-compute good initial guesses for the SNES solver when this fails to converge. This modification is necessary to solve \eqref{eq:NSweak} at high Reynolds-numbers. The code used to obtain these results can be found at~\cite{zenodo-firedrake-and-driver}.}

\begin{figure}[htb!]
    \centering
    \includegraphics[width=0.8\linewidth]{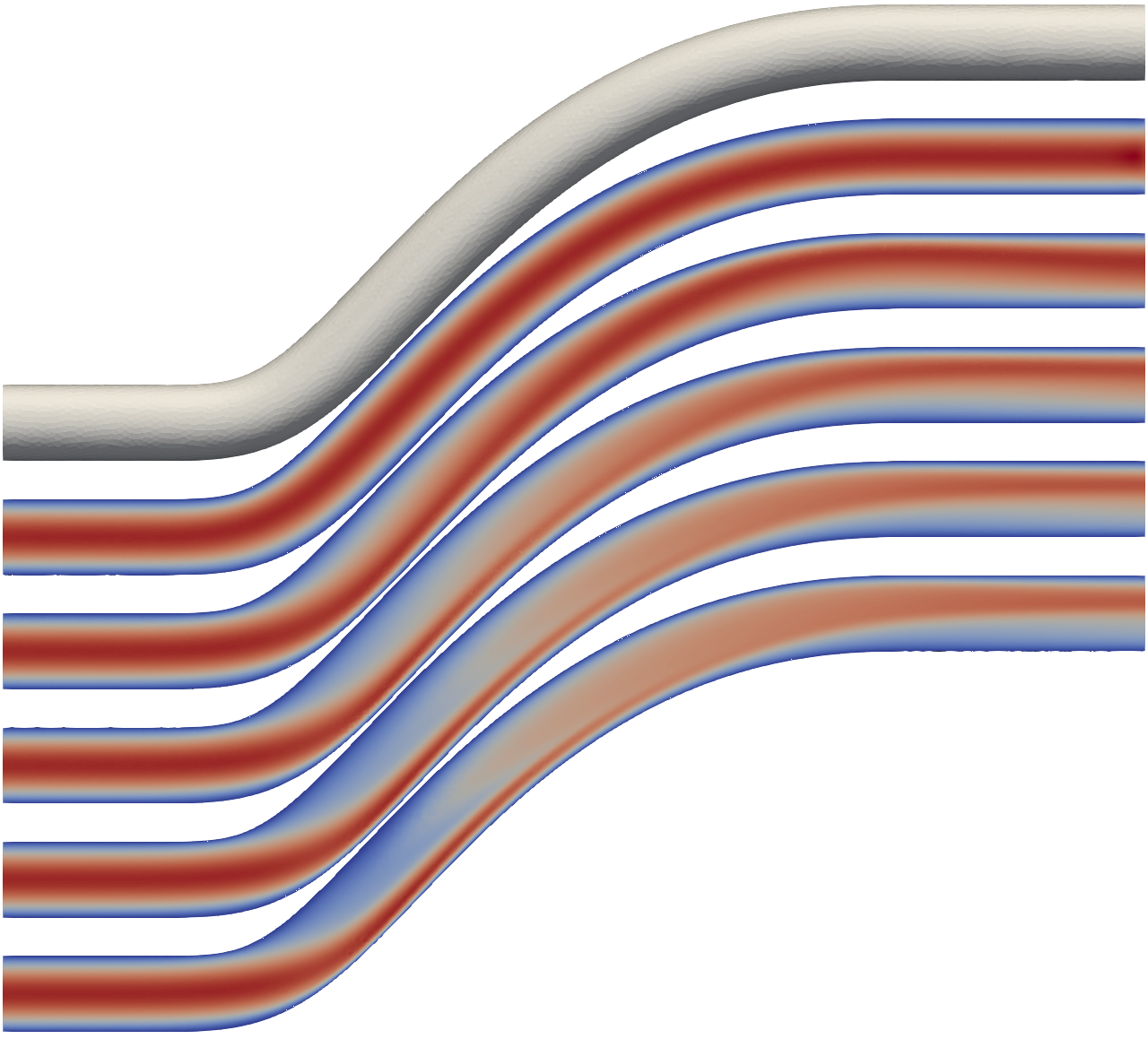}
\caption{Initial design of the pipe and velocity solutions for $\Re=1,250,500,750,1000$ (from top to bottom) inside the pipe. Red corresponds to high velocity, blue corresponds to low velocity.}
\label{fig:initial-shape}
\end{figure}

\begin{figure}[htb!]
    \centering
    \includegraphics[width=0.8\linewidth]{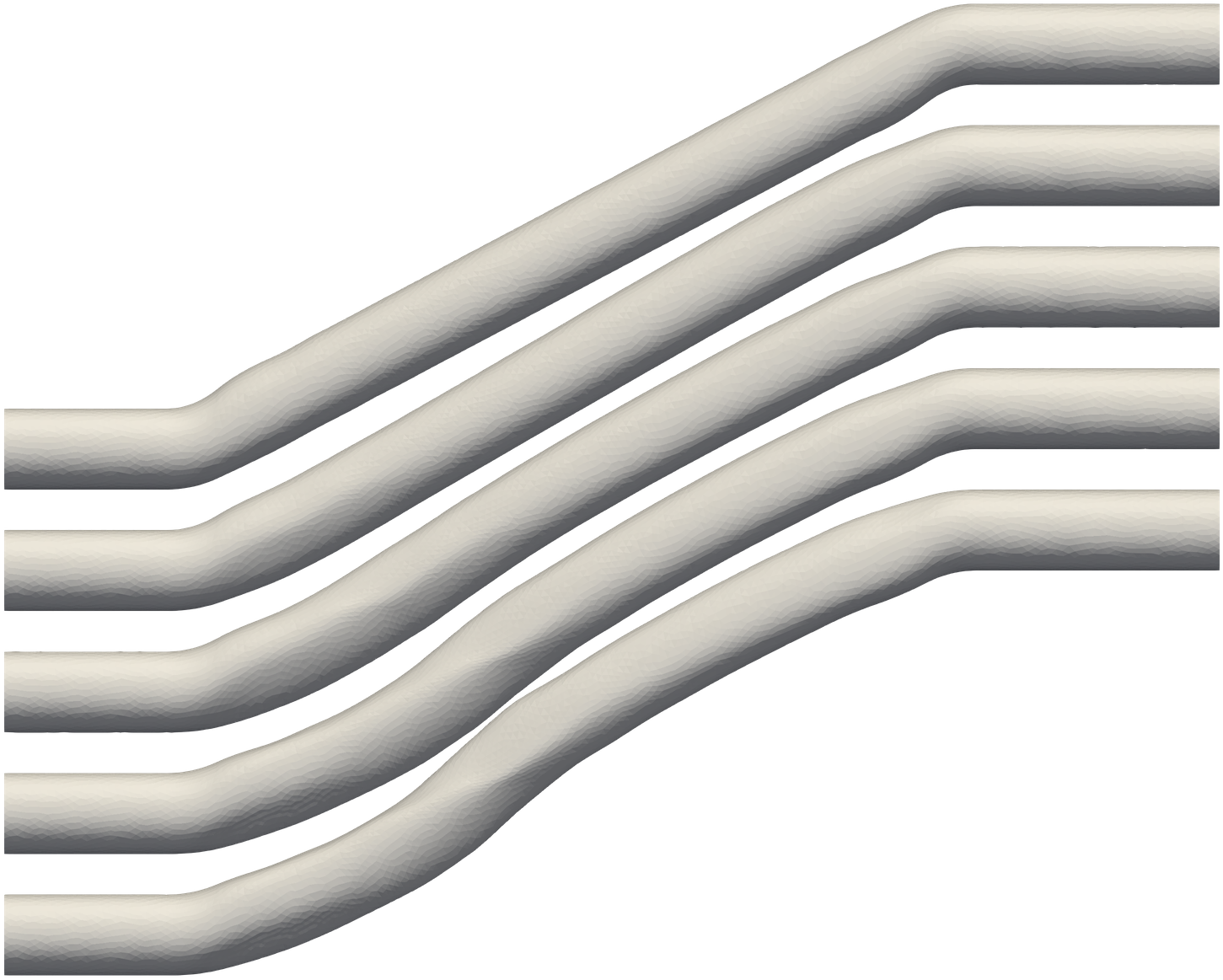}
    \caption{Optimized shapes for $\Re=1,250,500,750,1000$ (from top to bottom).}
\label{fig:optimal-shapes}
\end{figure}

\begin{figure}[htb!]
    \centering
    \includegraphics[width=0.8\linewidth]{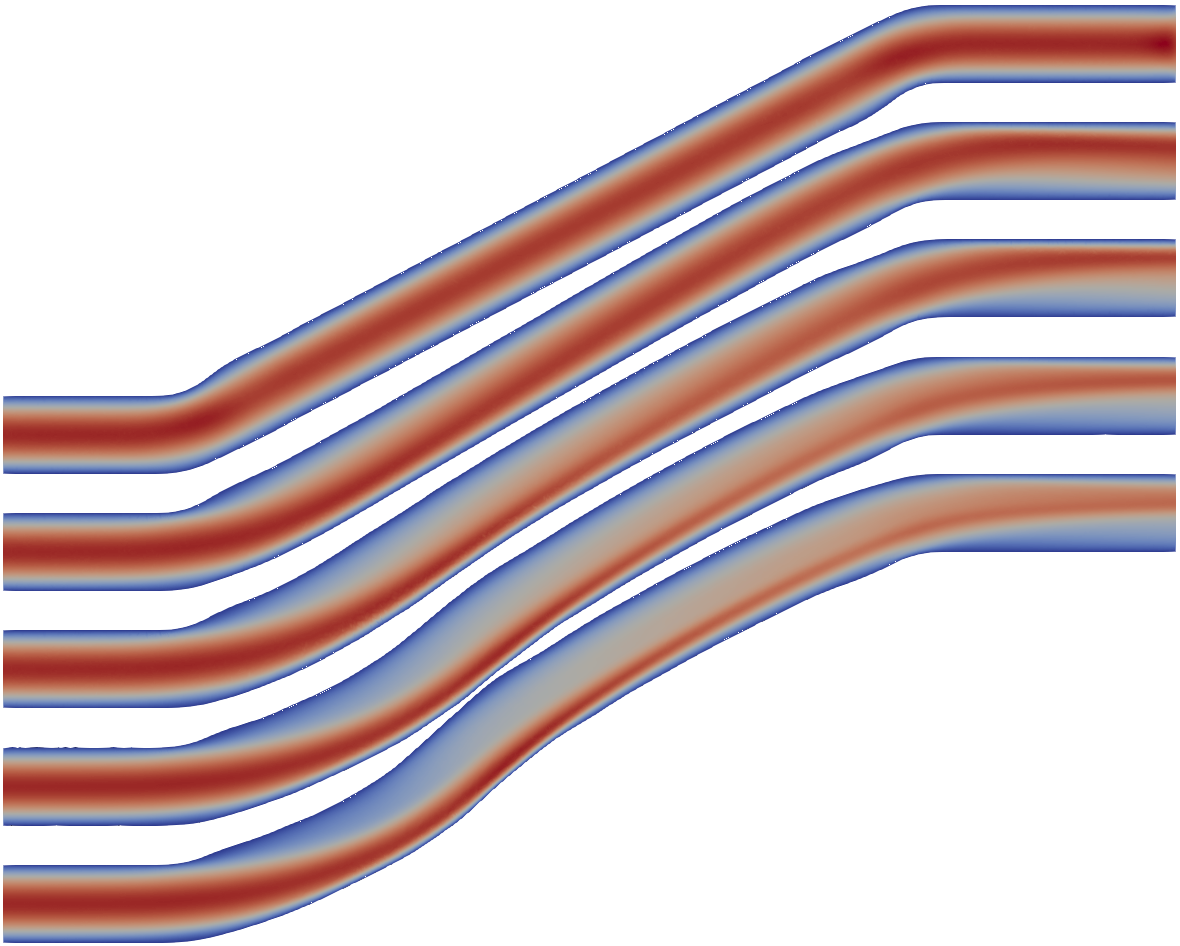}
    \caption{Velocity field for optimized shapes for $\Re=1,250,500,750,1000$ (from top to bottom). Red corresponds to high velocity, blue corresponds to low velocity.}
\label{fig:optimal-shapes-flow}
\end{figure}

In Figure~\ref{fig:initial-shape}, we show the initial design and plot the magnitude of the fluid velocity on a cross section of the pipe for different Reynold-numbers.
In Figure~\ref{fig:optimal-shapes}, we show the resulting optimized shapes, and in Figure~\ref{fig:optimal-shapes-flow} the corresponding magnitudes of the fluid velocity.
Qualitatively, we observe that, as the Reynolds-number increase, we obtain increasingly S-shaped designs that avoid high curvature at the two fixed ends.
Finally, we remark that the objective is reduced by approximately $6.9\%$, $9.5\%$, $10.8\%$, $11.7\%$ and $12.1\%$ \textcolor{black}{in at most
10 augmented Lagrangian steps}, respectively. \textcolor{black}{Each numerical
experiment was run in parallel on a server with 36 cores.}

{\color{black}
\section{Example: Shape optimization of a cantilever}
\label{sec:cantilever}

In this second tutorial, we consider another shape optimization classic
\cite{AlPa06, Al07}: the minimization of the compliance of a cantilever
$\rm{\Omega}$ (see Figure \ref{fig:cantilever-sketch}) subject to a dead
surface load $\Vg$. To make this example reproducible on a basic laptop
we consider a two dimensional version of this problem. Adapting the
code presented in this section to simulate a three dimensional cantilever is
straightforward.

\begin{figure}[htb!]
\center
\begin{tikzpicture}
    \begin{scope}[scale=0.3]
        \draw[black, fill=black!10] (0,-2)--(0,-4)--(0.9,-5)--(7.0,-6)--(20,-1)--(20,+1)--(10.7,+5.8)--(0.85,+5)--(0,+4)--(0,+2)--(0.5,+2.4)--(3.2,-0.2)--(0.5,-2.4)--(0,-2);
        \draw[black,fill=white] (2.8,+3.8)--(2.5,+2.3)--(6.5,-0.2)--(8.7,+3.8)--(2.8,+3.8);
        \draw[black,fill=white] (10.9,+4.1)--(7.0,-1.8)--(8.1,-3.5)--(17.1,-0.1)--(10.9,+4.1)--(10.9,+4.1);
        \node[right] at (-0.1, -3.3) {$\Gamma_1$};
        \node[right] at (-0.1, +3.3) {$\Gamma_1$};
        \node[left] at (20, 0) {$\Gamma_2$};
        \draw[black, very thick] (0, -4) -- (0, -2);
        \draw[black, very thick] (0, +4) -- (0, +2);
        \draw[->, red, ultra thick] (20, 1) -- (20, -1.2);
        \node[right, red] at (20, 0) {$\mathbf{g}$};
    \end{scope}
\end{tikzpicture}
\caption{A cantilever is attached to a wall on the left ($\Gamma_1$) and
pulled downwards on the right ($\Gamma_2$) by a dead surface load. A poor cantilever
design may fail to offer sufficient support.}

\label{fig:cantilever-sketch}
\end{figure}
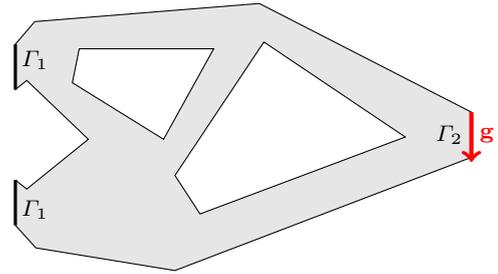

We model the cantilever elastic response using the linear elasticity
equations \cite[Ch. VI, Sect. 3]{Br07}. More precisely, we consider the
displacement formulation: the displacement field $\Vu$ satisfies the
variational problem

\begin{equation}\label{eq:linearElasticity}
\int_{\rm\Omega} \sigma(\Vu):\epsilon(\Vv) \DX
= \int_{\Gamma_2} \Vg\cdot\Vv \dS
\quad \text{ for every } \Vv\,, 
\end{equation}
where $\sigma(\Vu) = \lambda \Div(\Vu)\rm{I} + 2\mu\varepsilon(\Vu)$,
$\lambda$ and $\mu$ are the Lam\'{e} constants, $\rm{I}$ is the identity
matrix, $\varepsilon (\Vu) = \frac{1}{2}(\nabla \Vu+\nabla\Vu^\top)$, $\nabla
\Vu$ is the derivative (Jacobian matrix) of $\Vu$ and $\nabla\Vu^\top$ is the
derivative transposed, the colon symbol $:$ denotes the Frobenius inner
product, that is $\sigma(u):\varepsilon(\Vv)= trace(\sigma(\Vu)^\top
\varepsilon(\Vv))$, and $\Gamma_2$ denotes the part of the boundary that is
subject to the surface load $\Vg$. Additionally, although not explicitly
indicated in \eqref{eq:linearElasticity}, the displacement $\Vu$ vanishes on
$\Gamma_1$, which is the part of the boundary $\partial\rm\Omega$ that corresponds
to the wall, see Figure \ref{fig:cantilever-sketch}.

To model the compliance of the structure $\rm\Omega$,
we consider the functional

\begin{equation} \label{eq:compliance}
\rm{J}(\rm{\Omega}) =
\int_{\rm{\Omega}} \sigma (\Vu) : \varepsilon(\Vu) \DX\,,
\end{equation}
where $\Vu$ is the displacement and satisfies \eqref{eq:linearElasticity}.

As typical for this shape optimization test case, we also impose
a volume equality constraint on $\rm\Omega$ to prescribe the amount
of material used to fabricate the cantilever.

Similarly to Section \ref{sec:pipe}, in the following subsections
we explain how to solve this shape optimization test case in Fireshape.

\subsection{Step 1: Provide an initial guess}

To create the mesh of the initial guess, we use again the software Gmsh
\cite{GeRe09}. The initial guess employed is sketched in Figure
\ref{fig:cantilever-sketch}. The only mesh detail necessary to understand
this tutorial is that, in Listing \ref{lst:PDEconstraintCantilever}, the
number 1 corresponds to $\Gamma_1$, whereas the boundary flag 2 refers to
$\Gamma_2$. For the remaining geometric details, we refer to the code
archived on Zenodo \cite{zenodo-firedrake-and-driver}.

\subsection{Step 2: Implement the PDE-constraint}

In this example, the PDE constraint \eqref{eq:linearElasticity} is already in
weak form. To implement its finite element discretization, we create a class
\verb!LinearElasticitySolver! that inherits from Fireshape's class
\verb!PdeConstraint!, see Listing \ref{lst:PDEconstraintCantilever}. To
discretize \eqref{eq:linearElasticity}, we employ the standard piecewise
affine Lagrangian finite elements \cite[Ch. II, Sect. 5]{Br07}. To solve the
resulting linear system, we employ a direct solver. This is sufficiently fast
for the current 2D example.
For readers interested in creating a three dimensional
version of this example, we suggest replacing the direct solver with a Krylov method preconditioned by
one of the algebraic multigrid methods that can be accessed via PETSc.

\begin{mylisting}[htb!]
\includegraphics[width=\linewidth]{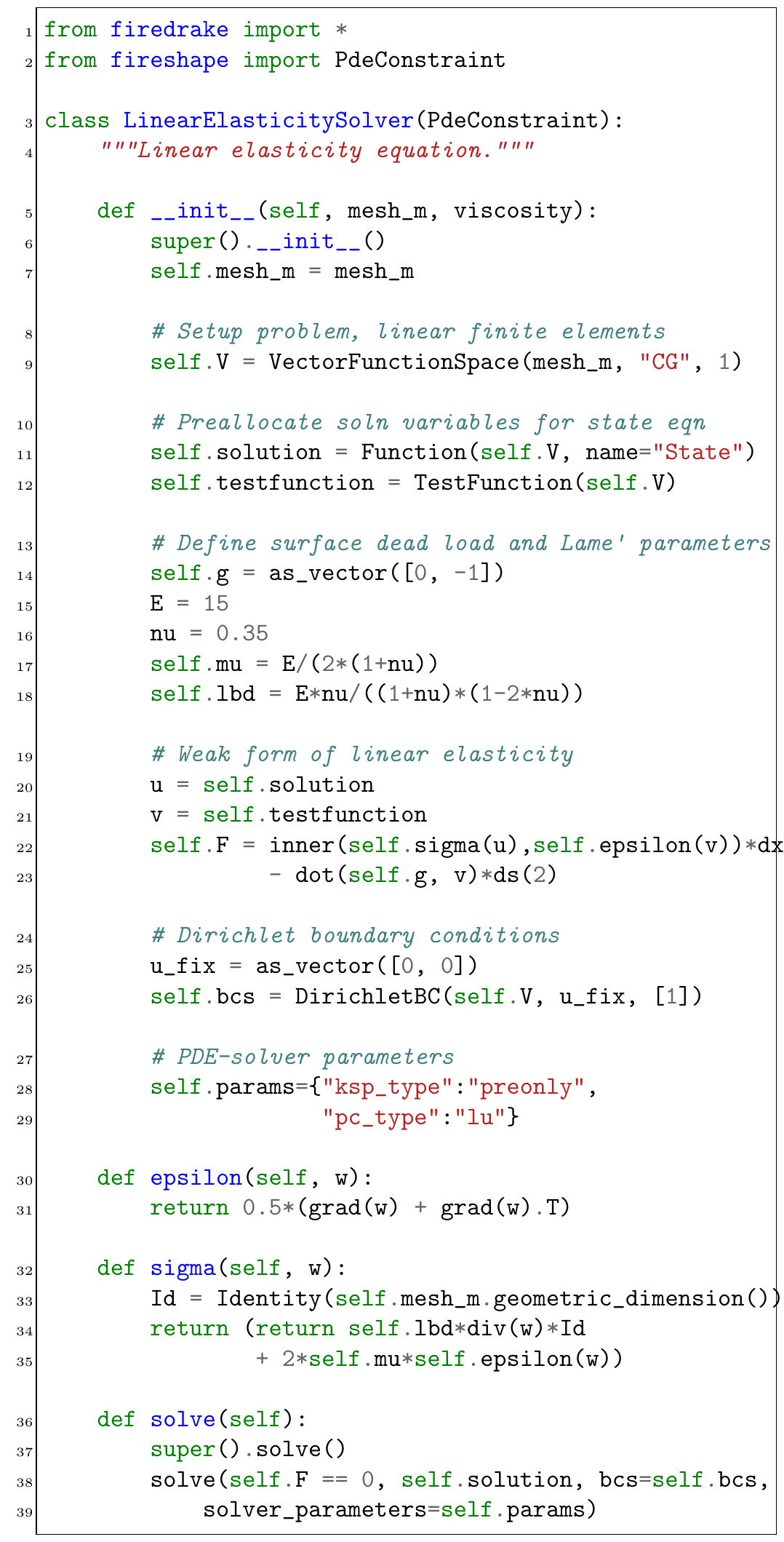}
\caption{Linear elasticity solver.}
\label{lst:PDEconstraintCantilever}
\end{mylisting}

\subsection{Step 3: Implement the objective function}

We implement the objective function $\rm{J}$ defined in Equation
\eqref{eq:compliance} in a class \texttt{Compliance} that inherits from
Fireshape's class \texttt{ShapeObjective}, see Listing
\ref{lst:objectiveCantilever}. This functional accesses the displacement
field $\Vu$ via the variable \texttt{self.pde\_solver}.

To verify that the domain used to compute $\Vu$ is feasible, that is, that
the mesh is not tangled, we create a piecewise constant function
\texttt{detDT} and interpolate the determinant of the derivative of the
transformation $\VT$ used to generated the updated domains (see Section
\ref{sec:shapeoptshapecalc} and \ref{sec:isoFEM} for more details about the
role of the transformation $\VT$ in Fireshape). If the minimum of
\texttt{detDT} is a small number, this indicates that the mesh has poor
quality (see also Remark \ref{rmk:detDTsmall}), and this negatively affects
the accuracy of the finite element method. In this case, we decide that the
value of $\rm{J}$ is \texttt{NaN} (``not a number'') to ensure that only
trustworthy simulations are used \footnote{In Section \ref{sec:step3}, we did
not include this extra test because the Navier-Stokes' solver usually fails
to converge on poor quality meshes.}.

\begin{mylisting}[htb!]
\includegraphics[width=\linewidth]{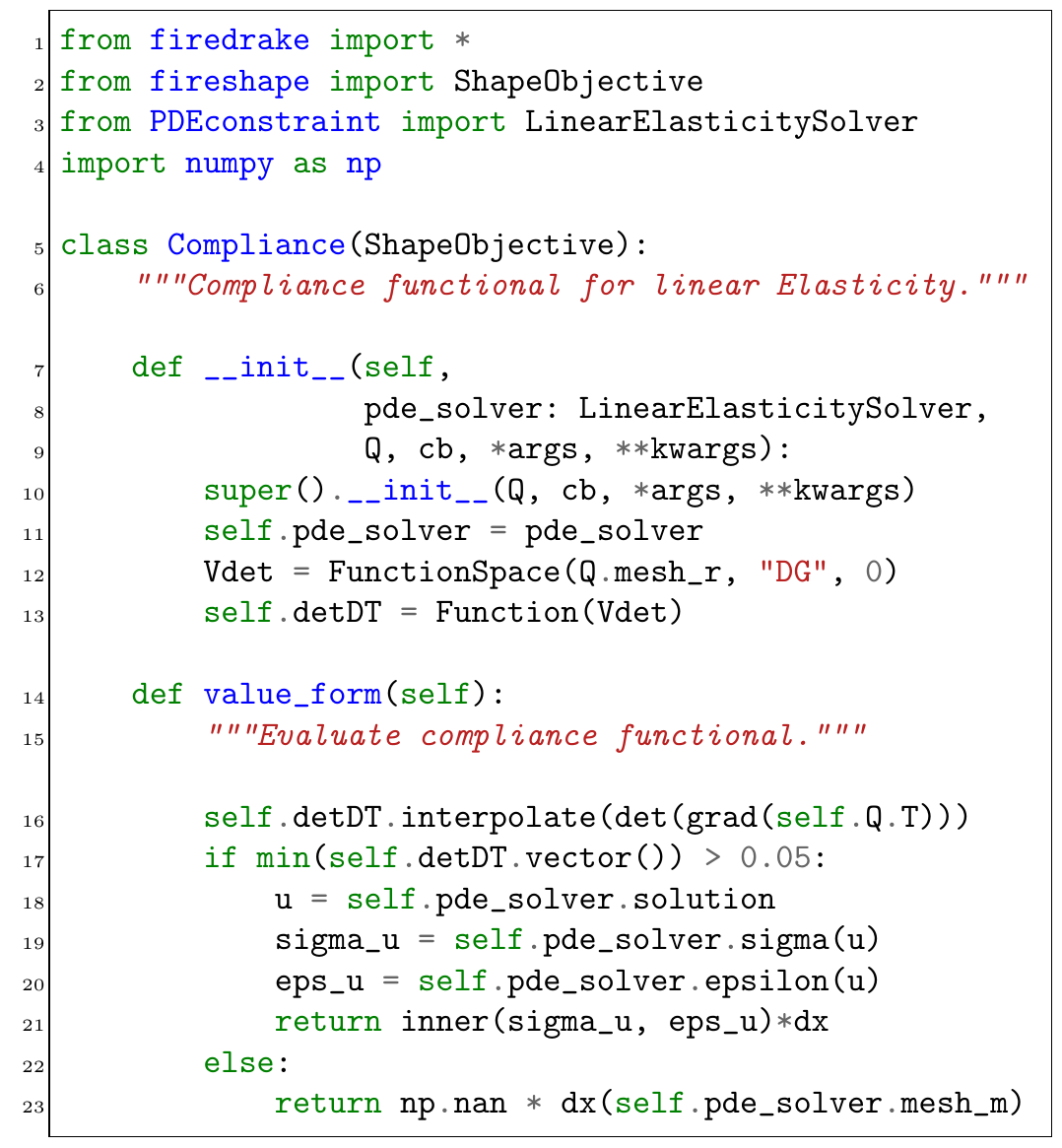}
\caption{Objective function $\rm J$ \eqref{eq:compliance},
which quantifies the compliance of the structure.
Note that $\rm J$ returns \texttt{NaN} (``not a number'')
if the mesh has poor quality.}
\label{lst:objectiveCantilever}
\end{mylisting}

\subsection{Final step: Set up and solve the problem}
\label{sec:finalstepCantilever}
Finally, to shape optimize the initial desing
we create a ``main file'' that loads the initial mesh, sets-up
the optimization problem, and solves it with an optimization algorithm.
The rest of this section contains a line-by-line description of
the ``main file'', which is listed in Listing \ref{lst:mainCantilever}.

\begin{mylisting}[htb!]
\includegraphics[width=\linewidth]{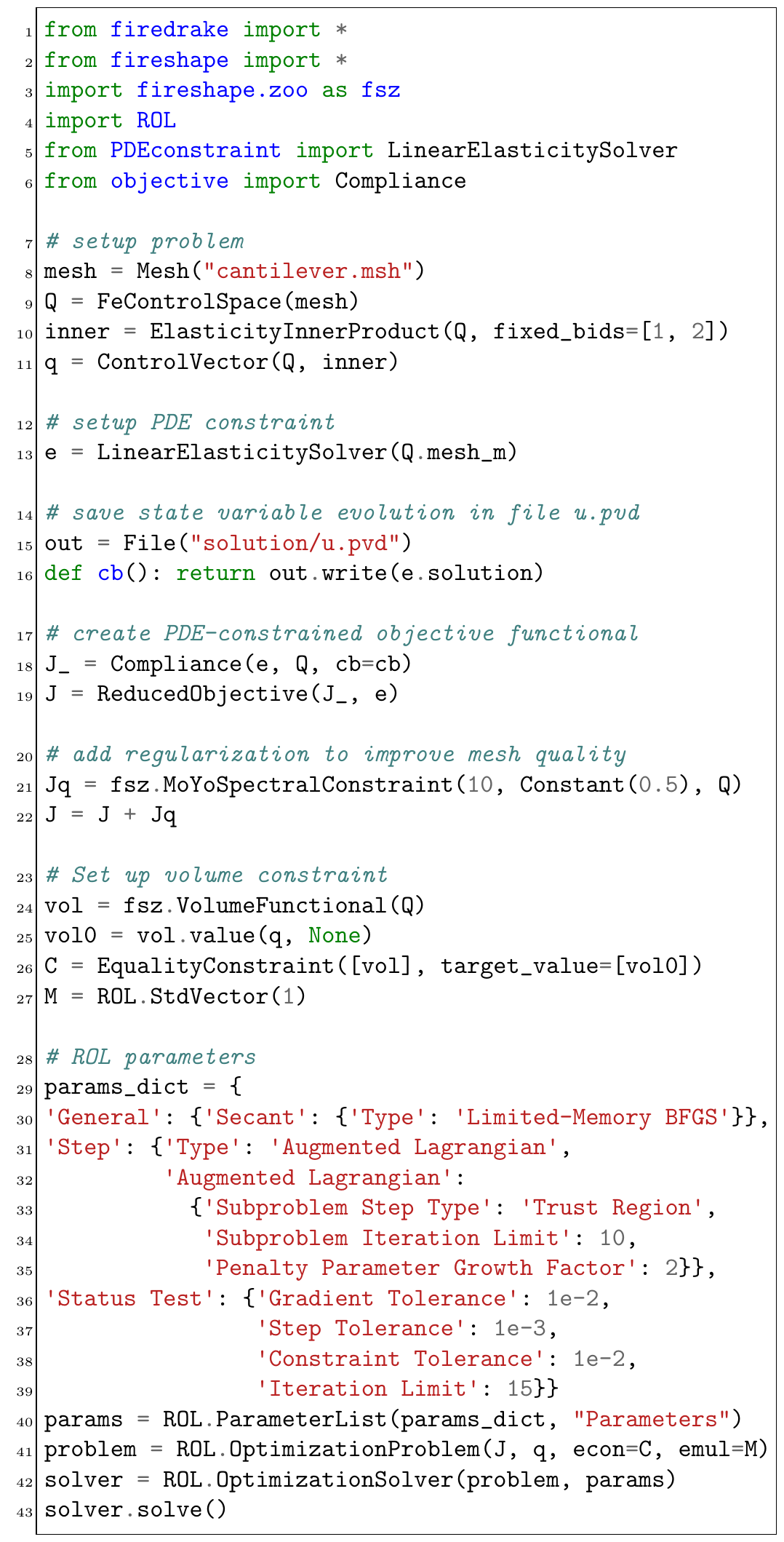}
\caption{Minimize the kinetic energy dissipation of a fluid using Fireshape.}
\label{lst:mainCantilever}
\end{mylisting}
\renewcommand\labelitemi{\textbullet}
\begin{itemize}
\item \emph{In lines 1-6,} we import all necessary Python libraries and modules.

\item \emph{In lines 8-11,} we load the mesh of the initial design and we
specify to discretize perturbations of the domain $\rm \Omega$ using
piecewise affine Lagrangian finite elements. We also employ a metric inspired
by the linear elasticity equations and specify that the boundaries $\Gamma_1$
and $\Gamma_2$ cannot be modified. Sections \ref{sec:controlspace} and
\ref{sec:innerproduct} \textcolor{black}{provide} more details about the
Fireshape classes \texttt{ControlSpace} and \texttt{InnerProduct}.

\item \emph{In line 13,} we initiate the linear elasticity solver.

\item \emph{In lines 15-16,} we tell Fireshape to store the finite element
approximation of the displacement field $\Vu$ in the folder \texttt{solution}
every time the domain $\rm \Omega$ is updated. We can visualize the evolution
of $\Vu$ along the optimization process by opening the file \texttt{u.pvd}
with Paraview \cite{ahrens2005paraview}.

\item \emph{In lines 18-19,} we initiate the objective function $\rm J$ and
the associated reduced functional, which is used by Fireshape to define the
appropriate Lagrange functionals to shape differentiate $\rm J$. We refer to
Section \ref{sec:shapeoptshapecalc} for more details about shape
differentiation using Lagrangians. Note that, as mentioned in Remark
\ref{rmk:shapediff}, Fireshape does not require users to derive (and
implement) shape \textcolor{black}{derivative} formulas by hand. The whole
shape differentiation process is automated using pyadjoint \cite{HaMiPaWe19,
DoMiFu20}.

\item \emph{In lines 21-22,} we add an additional regularization term to $\rm
J$ to promote mesh quality of domain updates. This regularization term
controls the pointwise maximum singular value of the geometric transformation
used to update $\rm \Omega$ \cite{wechsungthesis}. See Section
\ref{sec:shapeoptshapecalc}, Figure \ref{fig:Q}, and Remark
\ref{rmk:detDTsmall} for more information about the role of geometric
transformations in Fireshape.

\item \emph{In lines 24-27,} we set up an equality constraint
to ensure that the volume of the initial and the optimized domains
are equal.

\item \emph{In lines 29-39,} we select an optimization algorithm
from the optimization library ROL.
More specifically, we use an augmented Lagrangian approach
\cite[Ch. 17.3]{NoWr06} with limited-memory BFGS Hessian updates
\cite[Ch. 6.1]{NoWr06} to deal with the volume equality constraint,
and a Trust-Region algorithm \cite[Ch. 4.1]{NoWr06} to solve the intermediate
models generated by the augmented Lagrangian algorithm.
\item \emph{In lines 44-47,} finally, we gather all information and solve the problem.
\end{itemize}

\subsection{Results}

Figure \ref{fig:cantileverOptimized} shows the initial design (top) and the
design optimized using Listing \ref{lst:mainCantilever} (middle). At first
look, the result looks reasonable. However, after an inspection of ROL's
output, it becomes apparent that the optimization algorithm stops after 8
(augmented Lagrangian outer) steps because the optimization step lenght
becomes too small. This is due to a triangle on the top left corner of the
cantilever that becomes almost flat (see Figure
\ref{fig:cantileverMeshDetail}, left). This deterioration of the mesh quality
is detected by the implementation of the objective function (Listing
\ref{lst:objectiveCantilever}, line 21), which forces the trust-region
optimization algorithm to reduce the trust-region radius. After 8 iterations,
this leads to an optimization step lenght that is smaller than $10^{-4}$, and
the optimization algorithm terminates. This is a positive result: the
algorithm detects that something is wrong and stops the simulation. It is the
task of the user to understand what originates the problem and provide a
remedy. After a look at the iterates on Paraview, it becomes apparent that
optimization algorithms tries to move a node above the boundary $\Gamma_1$,
which is fixed, and this leads to the generation of an elongated triangle. A
difficulty with this part of the design was already experienced in
\cite{AlPa06}, where the authors suggest to ``arbitrarily set the shape
gradient to zero near the corners of the shape''. Here, we present an
alternative solution.

\begin{figure}[htb!]
\includegraphics[trim=200 360 100 330, clip, width=1\linewidth]
{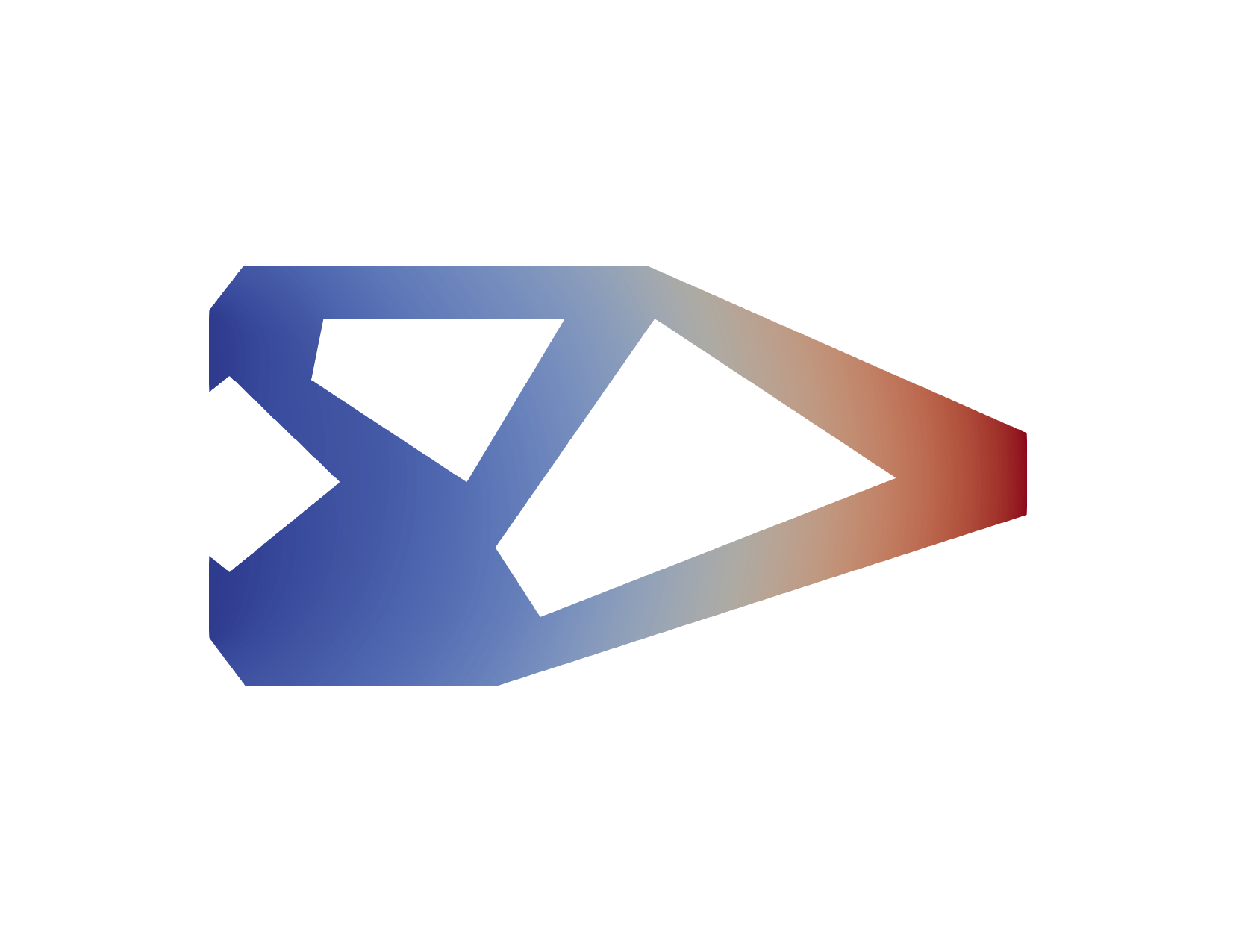}
\includegraphics[trim=200 330 100 330, clip, width=1\linewidth]
{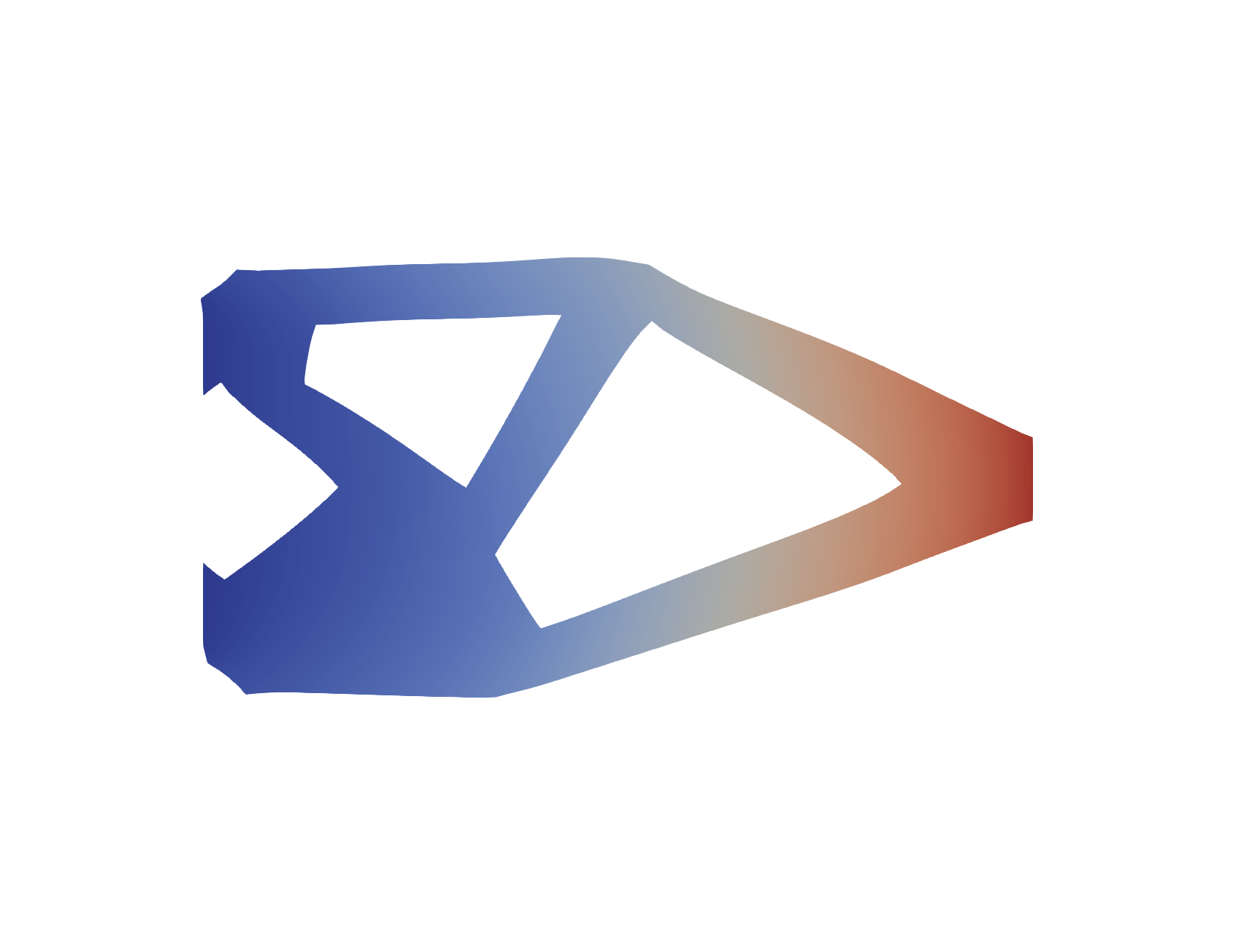}
\includegraphics[trim=200 330 100 330, clip, width=1\linewidth]
{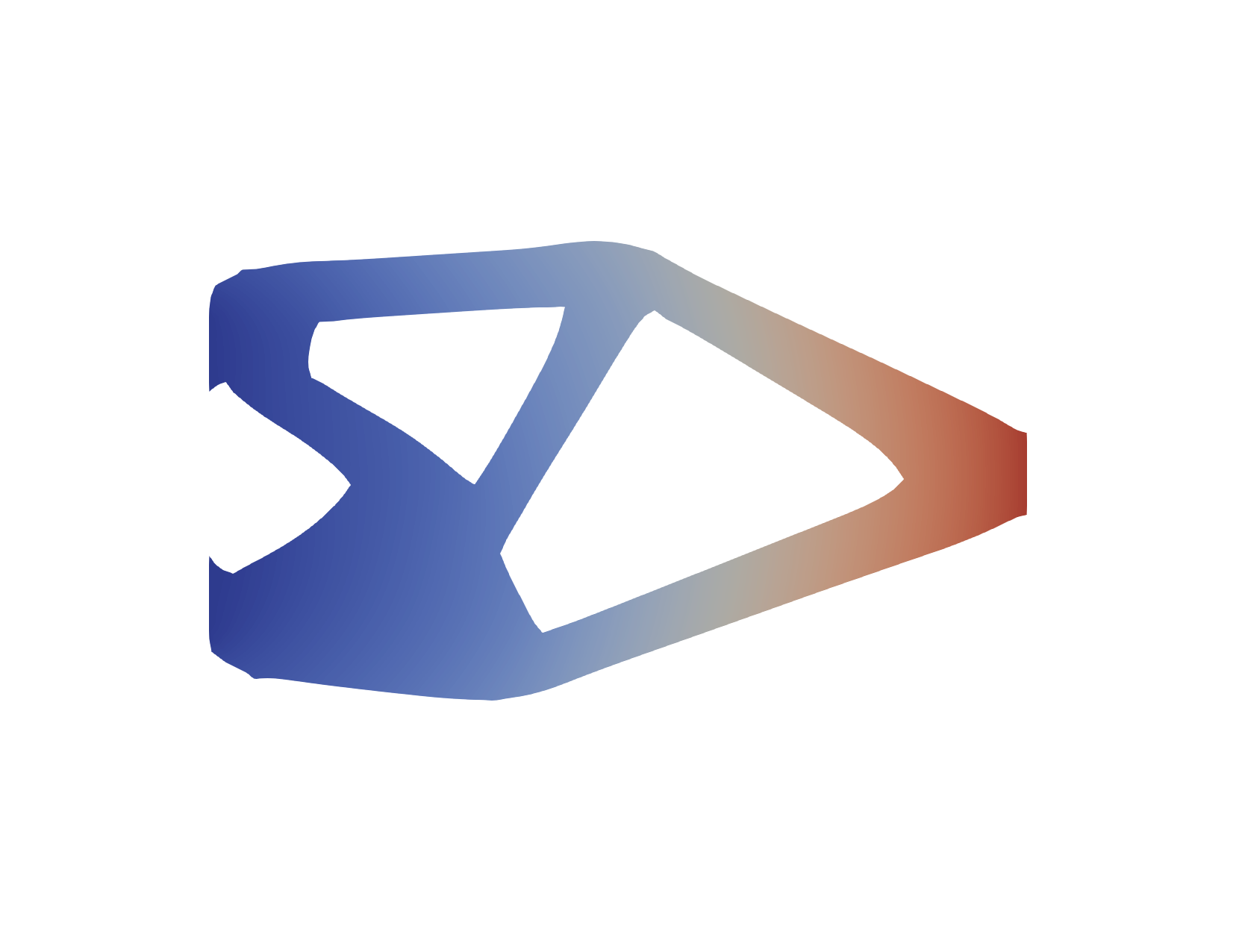}
\caption{\emph{Top:} Initial cantilever design. The absolute value of the displacement field $\Vu$ is colored from blue (minimum) to red (maximum).
\emph{Middle:} Cantilever designed optimized using Listing
\ref{lst:mainCantilever}. \emph{Bottom:} Cantilever designed optimized
by increasing the regularization of $\rm J$ in line 21 of Listing
\ref{lst:mainCantilever}.} 
\label{fig:cantileverOptimized}
\end{figure}

\begin{figure}[htb!]
\includegraphics[trim=300 300 800 200, clip, width=0.48\linewidth]
{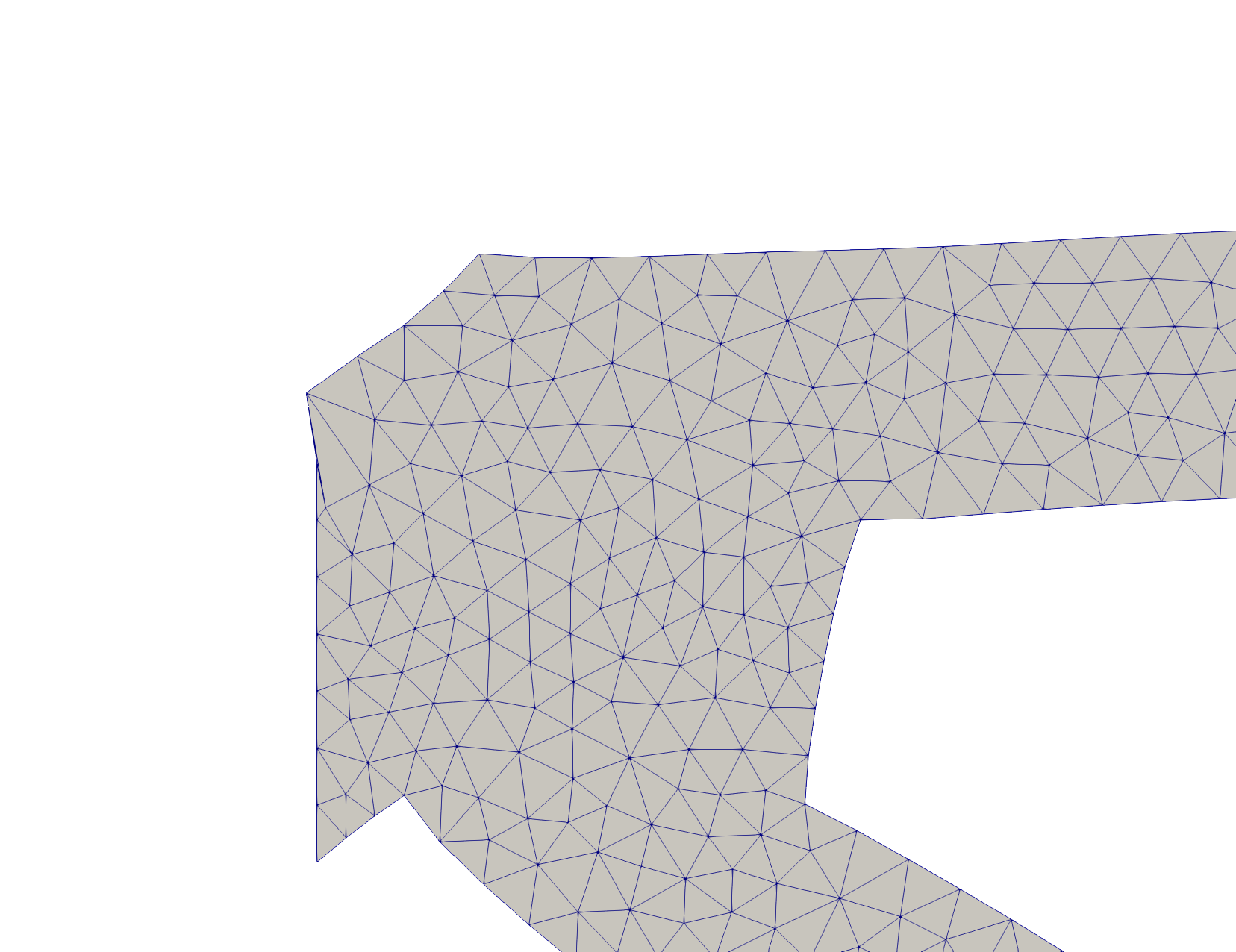}
\includegraphics[trim=300 300 800 200, clip, width=0.48\linewidth]
{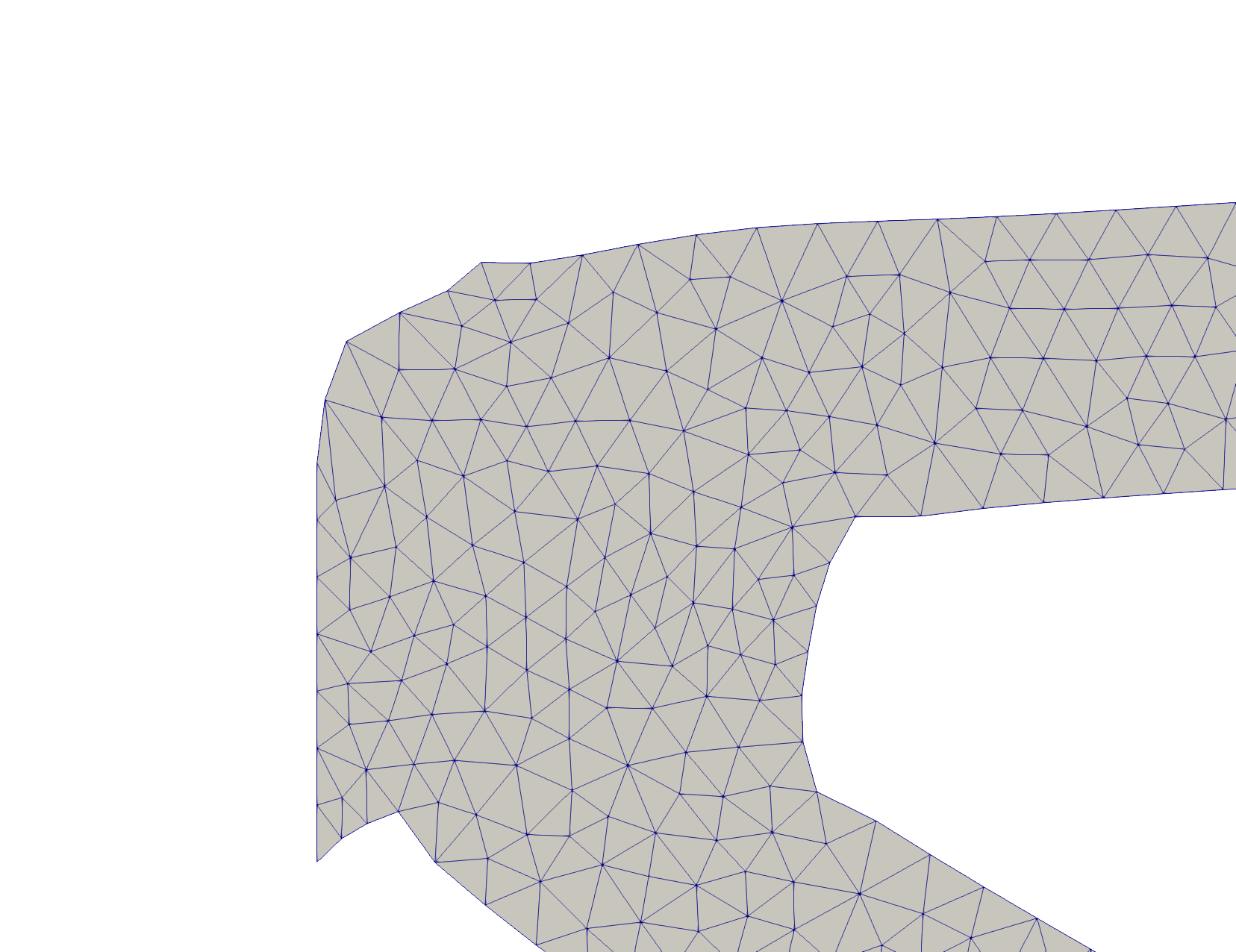}
\caption{\emph{Left:} Mesh detail of upper left corner of the cantilever
design optimized with Listing \ref{lst:mainCantilever}. The leftmost
triangle has almost collapsed to a line. \emph{Right:} Same mesh detail
when the regularization of $\rm J$ in line 21 of Listing
\ref{lst:mainCantilever} is increased. The elongated triangle is no
longer present.} 
\label{fig:cantileverMeshDetail}
\end{figure}

As mentioned in the previous paragraph, the issue is that the shape updates lead to the
creation of an elongated triangle. To overcome this problem, we can simply
increase from 10 to 100 the coefficient of the regularization term in line 21
Listing \ref{lst:mainCantilever}. The increased regularization sufficiently
penalizes this unwanted behavior, and the algorithm can further optimize the
design (see Figure \ref{fig:cantileverConvergence}).

\begin{figure}[htb!]
\includegraphics{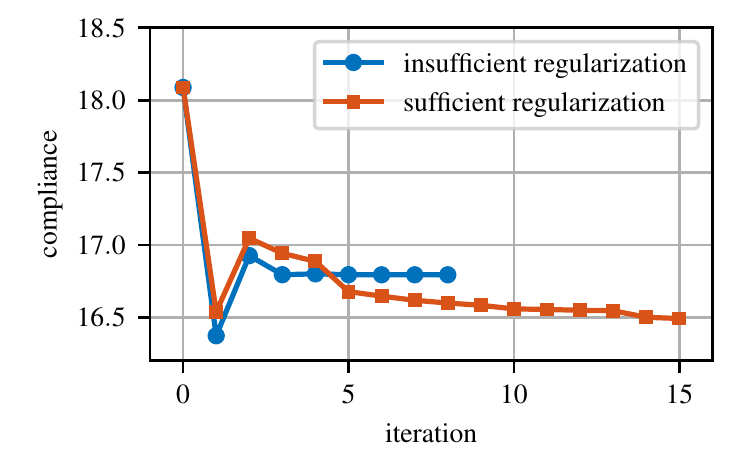}
\caption{Convergence history of the compliance. When the regularization
of $\rm J$ in line 21 of Listing \ref{lst:mainCantilever} is insufficient,
the algorithm stops after eight iterations because a triangle becomes very
elongated. However, if sufficient regularization is provided, the algorithm
does not experience this issue and further optimizes the design. Finally, we
point out that the small value of the compliance after the first iteration is
obtained by excessively violating the volume equality constraint. This is
subsequently corrected by the augmented Lagrangian algorithm.}
\label{fig:cantileverConvergence}
\end{figure}

}

\section{Shape optimization via diffeomorphisms}
\label{sec:shapeoptimization}

In this section, we describe the theory that underpins Fireshape.
We begin with a brief introduction to PDE-constrained optimization
to set the notation and mention the main idea behind optimization algorithms
(Subsection \ref{sec:PDEoptimization}). Then, we continue with an introduction
to shape calculus (Subsection \ref{sec:shapeoptshapecalc}).
Finally, we conclude with a discussion about the link between shape optimization and
parametric finite elements (Subsection \ref{sec:isoFEM}).

\subsection{Optimization with PDE-constraints}
\label{sec:PDEoptimization}

The basic ingredients to formulate a PDE-constrained optimization
\textcolor{black}{problem} are: a control variable $\Vq$ that lives in a
control space $\VQ$, a state variable $\Vu$ that lives in a state space $\VU$
and that solves a (possibly nonlinear) state equation $\VA(\Vq,\Vu) =
\mathbf{0}$, and a real function $\Jmat:\VQ\times\VU\to\bbR$ to be minimized.

\textbf{Example:}
\emph{Equation \eqref{eq:shapeoptproblem} is a PDE-constrained
optimization problem. The control variable $\Vq$ corresponds to 
the domain $\rm \Omega$, the state variable $\Vu$ represents the pair
velocity-pressure $(\Vu,\rm p)$, the nonlinear constraint $\VA$ represents
the Navier-Stokes equations \eqref{eq:NavierStokes}, and $\rm J$ denotes
the objective function \eqref{eq:kineticfunctional}. The state space
$\VU$ corresponds to the space of pairs $(\Va,\rm b)$ with weakly differentiable
velocities $\Va$ that satisfy $\Va = \Vg$ on $\partial \rm \Omega \setminus \Gamma$
and square integrable pressures $\rm b$ \cite[Ch. 8.2]{ElSiWa14}.
The control space $\VQ$ is specified in Section \ref{sec:shapeoptshapecalc}; see
Equation \eqref{eq:Q}.}

Most\footnote{An alternative approach is to employ so-called
``one-shot'' methods, where the optimality system of the problem
is solved directly \cite{Sc09}. Note that, since the optimality system
is often nonlinear, one-shot methods still involve iterative algorithms.}
numerical methods for PDE-constrained optimization
attempt to construct a sequence of controls
$\{\Vq^{(k)}\}_{k\in\bbN}$ and corresponding states
$\{\Vu^{(k)}\}_{k\in\bbN}$ such that
\begin{equation*}
\lim_{k\to\infty} \Jmat(\Vq^{(k)},\Vu^{(k)}) =
\inf_{\Vq,\Vu} \{\Jmat(\Vq,\Vu) : \VA(\Vq,\Vu) = \mathbf{0}\}\,.
\end{equation*}
Often, the sequence $\{\Vq^{(k)}\}_{k\in\bbN}$
is constructed using derivatives of the function $\Jmat$
and of the constraint $\VA$.
Common approaches are steepest descent algorithms
and Newton methods (in their quasi, Krylov, or
semi-smooth versions) \cite{HiPiUlUl09,Ke03}.
These algorithms ensure that the sequence
$\{\Jmat(\Vq^{(k)},\Vu^{(k)})\}_{k\in\bbN}$ decreases.
In special cases, it is even possible to show that the sequence
$\{\Vq^{(k)}\}_{k\in\bbN}$ converges \cite{HiPiUlUl09}. Although
it may be difficult to ensure these assumptions are met in industrial applications,
these optimization algorithms are still powerful tools to improve
state-of-the-art designs and are widely used to perform shape optimization.

\subsection{Shape optimization and shape calculus}
\label{sec:shapeoptshapecalc}
Shape optimization with PDE constraints is a particular branch of
PDE-constrained optimization where the control space $\VQ$
is a set of domain shapes. The space of shapes is notoriously difficult
to characterize uniquely. For instance, one could describe shapes through their
boundary regularity, or as level-sets, or as local epigraphs \cite[ch. 2]{DeZo11}.
The choice of the shapes' space characterization plays an important
role in the concrete implementation of a shape optimization algorithm,
and it can also affect formulas that result by differentiating
$\Jmat$ with respect to perturbations of the shape $\Vq$. To be more precise,
different methods generally lead to the same first order shape derivative formula,
but differ on shape derivatives of higher order \cite[ch. 9]{DeZo11}

\begin{figure}[htb!]
\centering
\begin{tikzpicture}
\node at (0,0) {\includegraphics[width = 1.2cm]{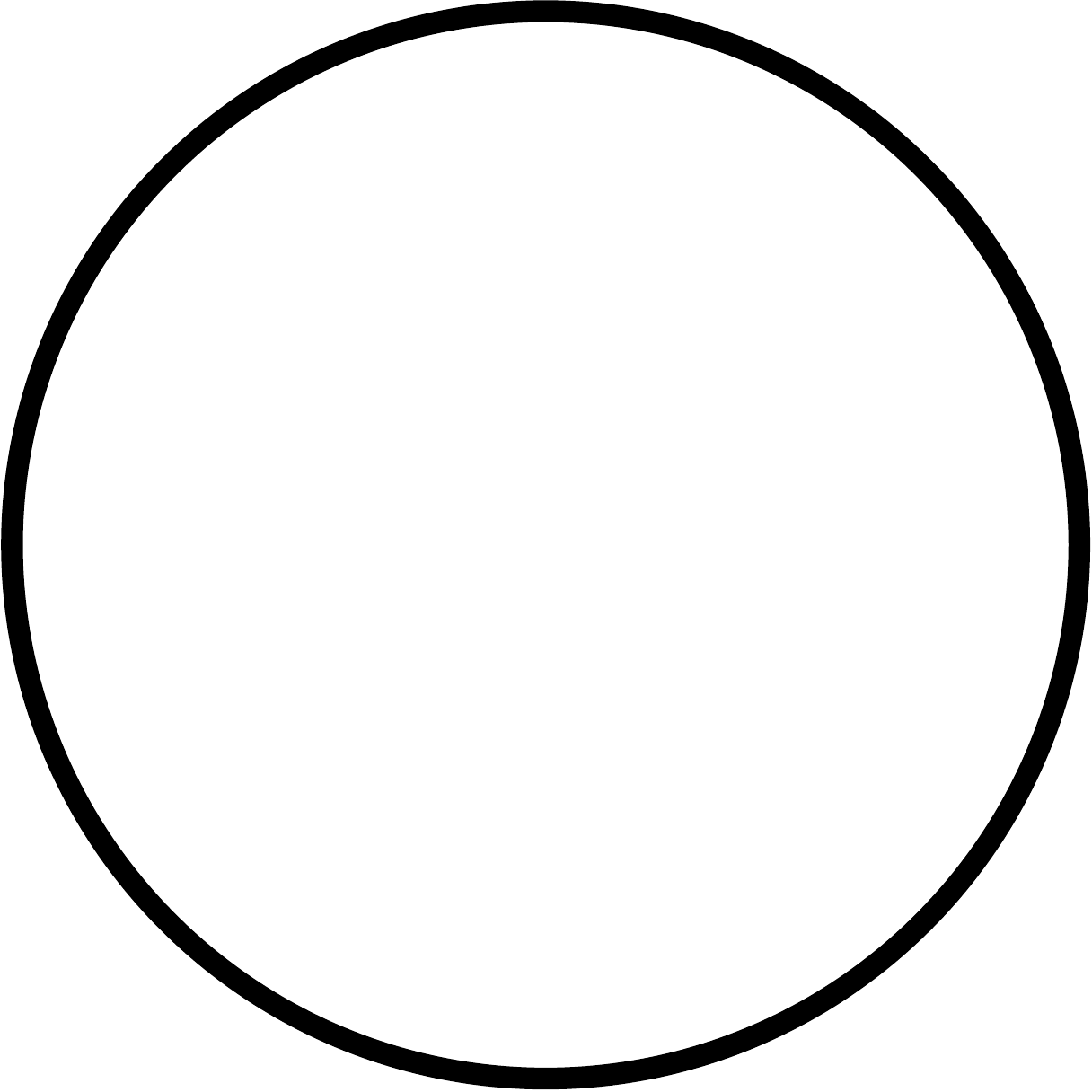}};
\node at (25:3) {\includegraphics[width = 1cm]{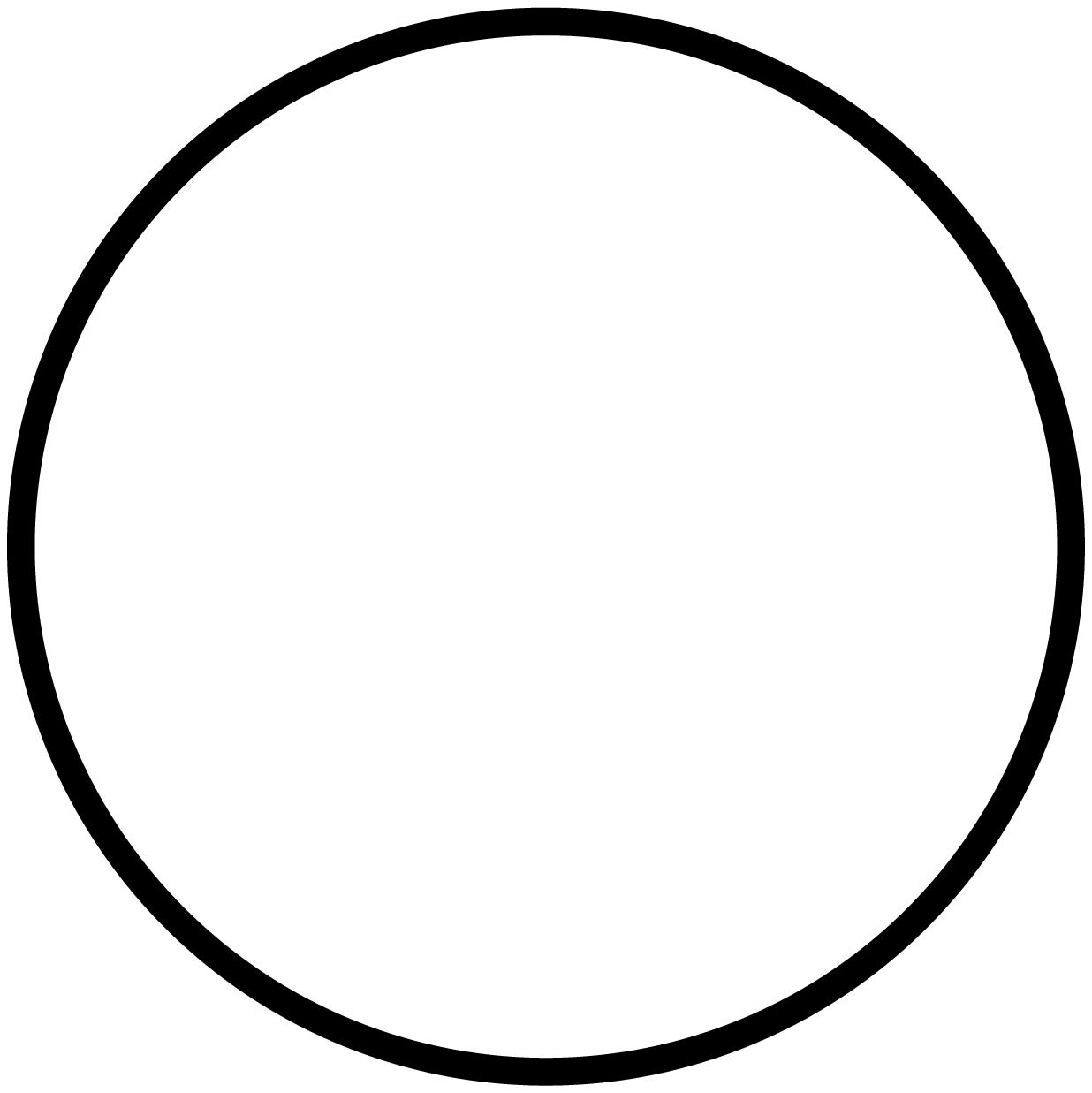}};
\node at (3,0) {\includegraphics[width = 1.6cm]{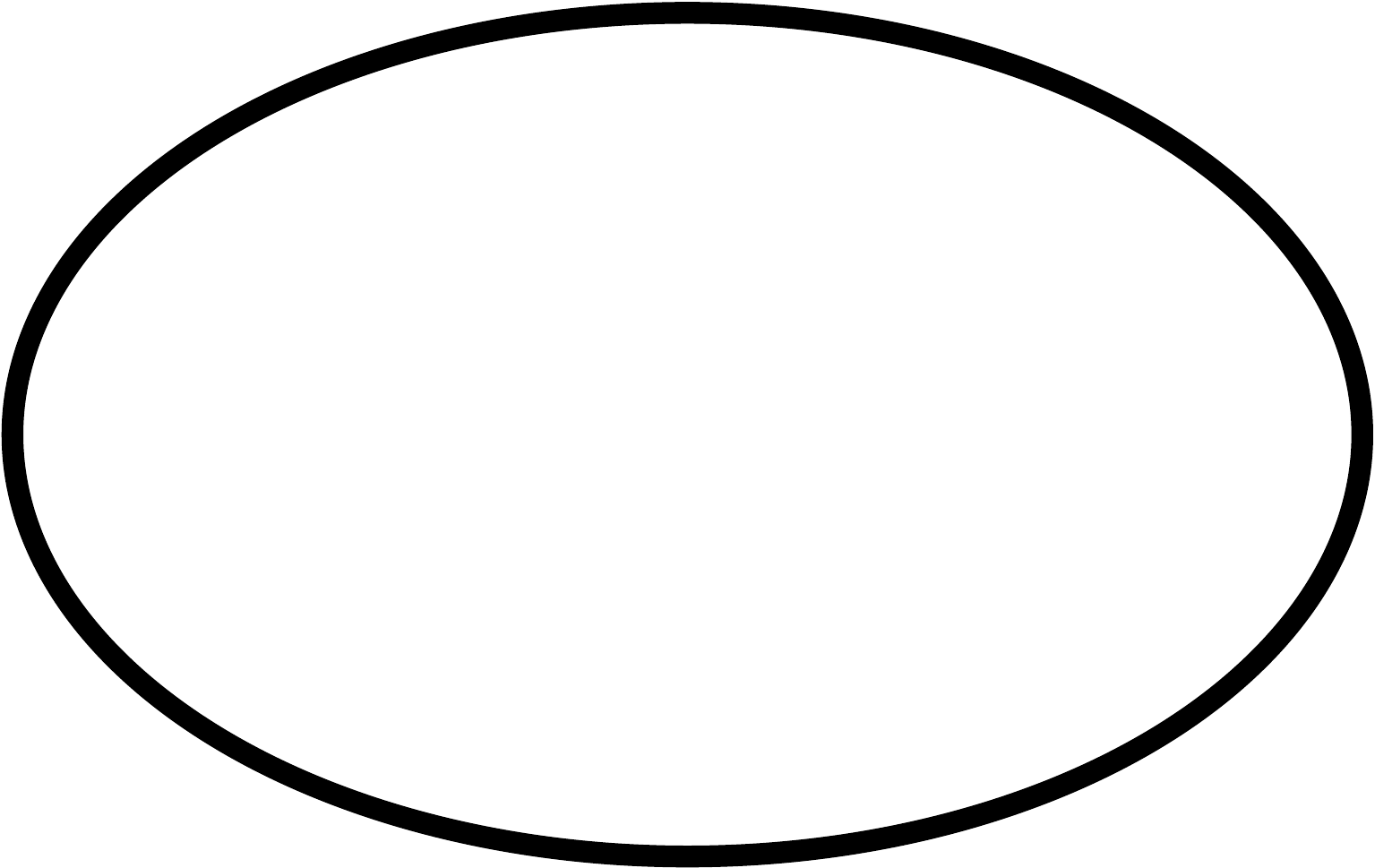}};
\node at (-25:3) {\includegraphics[width = 1.2cm]{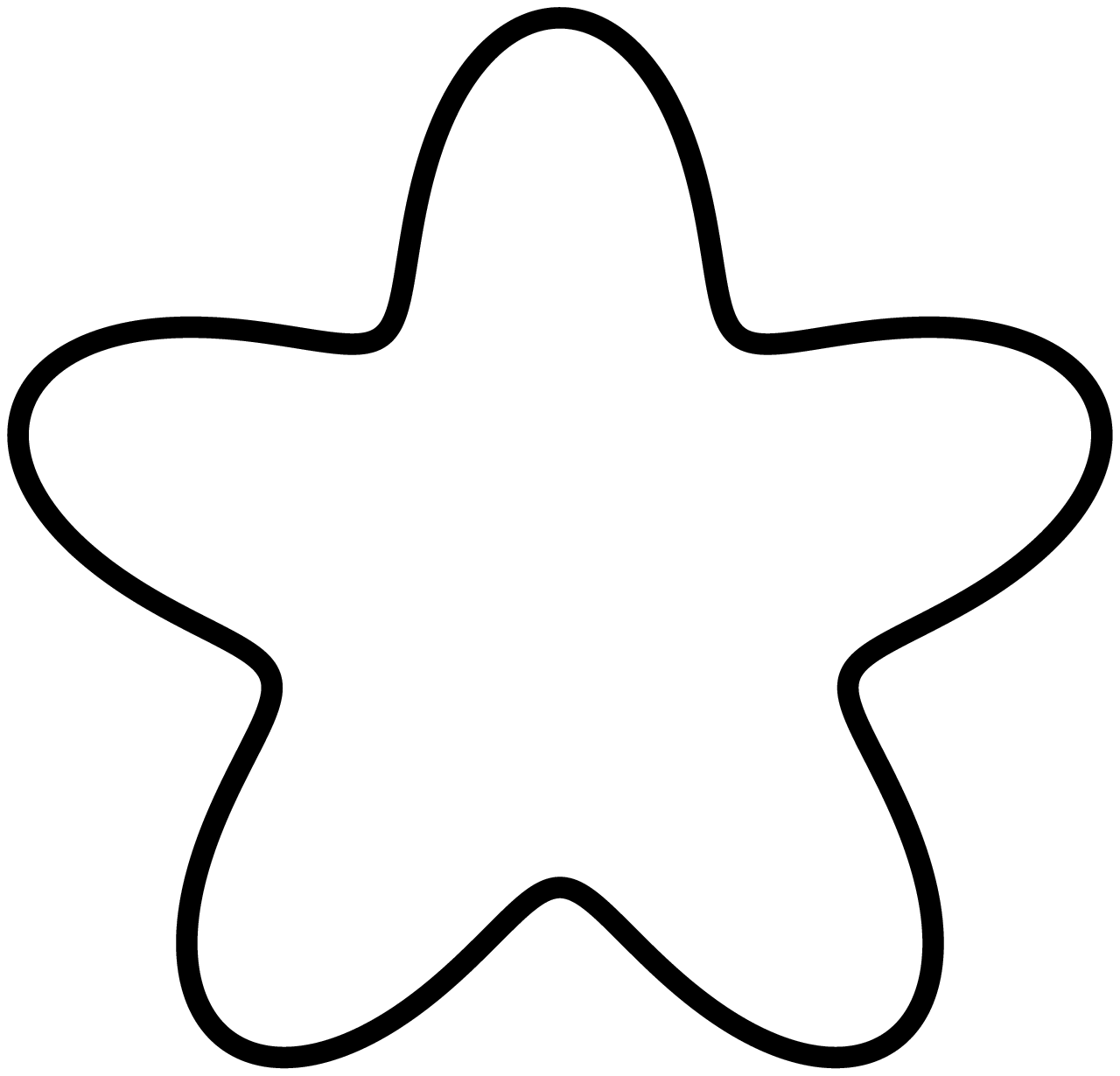}};
\node at (0,0) {$\rm{\Omega}$};
\draw[->, very thick] (25:1) -- (25:2);
\draw[->, very thick] (1,0) -- (2,0);
\draw[->, very thick] (-25:1) -- (-25:2);
\draw[->, very thick] (-50:1) -- (-50:2);
\node at (-50:2.3) {...};
);
\end{tikzpicture}
\caption{To construct the control space $\VQ$, we select an initial set
$\rm{\Omega}$ and a set of diffeomorphisms $\{\VT\}$ (\emph{black arrows}),
and collect every image $\VT(\rm{\Omega})$.}
\label{fig:Q}
\end{figure}

In this work, we model the control space $\VQ$ as the images of
bi-Lipschitz geometric transformations $\VT$ applied to an initial set
$\rm{\Omega}\subset\bbR^d$ \cite[ch. 3]{DeZo11}, that is,
\begin{equation}\label{eq:Q}
\VQ \coloneqq
\{\Vq = \VT(\rm{\Omega}) : \VT:\bbR^d\to\bbR^d \text{ is bi-Lipschitz} \}\,,
\end{equation}
see Figure \ref{fig:Q}.
We choose this model of $\VQ$ because it provides an explicit
description of the domain boundaries via
$\partial(\VT(\rm{\Omega})) = \VT(\partial\rm{\Omega})$,
and because it is compatible with higher-order finite elements \cite{PaWeFa18},
as explained in detail in the Section \ref{sec:isoFEM}.
Henceforth, we use the term differomorphism to indicate that a geometric transformation $\VT$ is bi-Lipschitz.

In this setting, the shape derivative of $\Jmat$ corresponds to the classical
G\^{a}teaux derivative in the Sobolev space
$W^{1,\infty}(\bbR^d,\bbR^d)$ \cite[Def. 1.29]{HiPiUlUl09}.
To see this, let us momentarily remove the PDE-constraint, and only consider a shape
functional of the form $\rm{\Omega}\mapsto\Jmat(\rm{\Omega})$.
The shape derivative of $\Jmat$ at $\rm{\Omega}$ is the
linear and continuous operator
$\rm{dJ}(\Omega, \cdot):W^{1,\infty}(\bbR^d,\bbR^d)\to\bbR$ defined by
\begin{equation}\label{eq:shapeder}
\rm{dJ}(\Omega,\VV)
\coloneqq \lim_{t\searrow 0}\frac{\Jmat(\Omega_t) - \Jmat(\Omega)}{t}\,,
\end{equation}
where $\rm{\Omega_t}$ is the set
$\rm{\Omega_t} = \{\Vx + \rm{t}\VV(\Vx) : \Vx\in \rm{\Omega}\}$
\cite[Def. 6.15]{Al07}.
By replacing $\rm \Omega$ with $\Vq = \VI(\rm \Omega)$, where $\VI$ denotes
the identity transformation defined by $\VI(\Vx) = \Vx$ for any
$\Vx\in \bbR^d$, Equation \eqref{eq:shapeder} can be equivalently rewritten as
\begin{equation*}
\rm{dJ}(\Vq,\VV)
\coloneqq \lim_{t\searrow 0}\frac{\Jmat(\Vq+t\VV) - \Jmat(\Vq)}{t}\,.
\end{equation*}
We highlight that this interpretation immediately generalizes to any $\Vq$ in $\VQ$
and can be used to define higher order shape derivatives.

The same definition \eqref{eq:shapeder}
of shape derivative holds in the presence of PDE-constraints:
the shape derivative at $\Vq$ of the function $\Jmat:\VQ\times\VU\to\bbR$ subject to
$\VA(\Vq,\Vu) = \mathbf{0}$ is the linear and continuous operator defined by
\begin{equation}\label{eq:shapederPDE}
\rm{dJ}(\Vq,\Vu,\VV)
\coloneqq \lim_{t\searrow 0}\frac{\Jmat(\Vq+t\VV,\Vu_t) - \Jmat(\Vq,\Vu)}{t}\,,
\end{equation}
where $\Vu_{\rm t}$ is the solution to $\VA(\Vq+\rm{t}\VV,\Vu_{\rm t})=\mathbf{0}$.
Computing shape derivative formulas using \eqref{eq:shapederPDE}
may present the difficulty of computing the shape derivative of $\Vu$, which
intuitively arises by the ``chain rule'' formula. The shape derivative of $\Vu$,
which is sometimes called ``material derivative'' of $\Vu$, can be eliminated
with the help of adjoint equations \cite[Sec. 1.6.2]{HiPiUlUl09}.
This process can be automated by introducing
the Lagrange functional \cite[Sec. 1.6.3]{HiPiUlUl09}
\begin{equation}\label{eq:Lagrangian}
\Cl(\Vq,\Vu,\Vp) \coloneqq \Jmat(\Vq,\Vu) + \langle \VA(\Vq,\Vu),\Vp\rangle\,.
\end{equation}
The term $\langle \VA(\Vq,\Vu),\Vp\rangle$ stems from testing the equation
$\VA(\Vq,\Vu) = \mathbf{0}$ with a test function $\Vp$ in the same way
it is usually done when writing a PDE in its weak form.

\textbf{Example: }
\emph{If $\VA$ denotes the Navier-Stokes equations \eqref{eq:NavierStokes},
then $\langle \VA(\Vq,\Vu),\Vp\rangle$ corresponds to the weak formulation \eqref{eq:NSweak}.}

The advantage of introducing the Lagrangian \eqref{eq:Lagrangian} is that,
by choosing $\Vp$ as
the solution to the adjoint equation
\begin{equation*}
\langle \partial_\Vu\VA(\Vq,\Vu),\Vp\rangle = -\partial_\Vu \Jmat(\Vq,\Vu)\quad
\text{for all } \Vu \text{ in } \VU\,,
\end{equation*}
where $\partial_\Vu$ denotes the partial differentiation
with respect to the variable $\Vu$, the shape derivative of $\Jmat$
can be computed as
\begin{equation*}
\rm{dJ}(\Vq,\Vu,\VV) = \partial_\Vq\Cl(\Vq,\Vu,\Vp,\VV)\,,
\end{equation*}
where $\partial_\Vq$ denotes partial differentiation
with respect to the variable $\Vq$. This is advantageous because
partial differentiation does not require computing the shape derivative of $\Vu$.

The Lagrangian approach to compute derivatives of PDE-constrained functionals
is well established \cite{ItKu08}, and its steps can be replicated by
symbolic computation software. Probably, the biggest success in this
direction is the pyadjoint project \cite{FaHaFuRo13,MiFuDo19}, which
derives ``the adjoint and tangent-linear equations and solves them using the
existing solver methods in FEniCS/Firedrake" \cite{MiFuDo19}. Thanks to the
shape differentiation capabilities of UFL introduced in \cite{HaMiPaWe19},
pyadjoint is also capable of shape differentiating PDE-constrained
functionals \cite{DoMiFu20}.

\begin{remark}\label{rmk:metric}
Optimization algorithms are usually based on 
steepest descent directions to update the control variable $\Vq$.
A steepest descent direction is a direction $\VV^*$ that minimizes the
derivative $\rm{dJ}$. Since $\rm{dJ}$ is linear, it is necessary to introduce
a norm $\Vert\cdot\Vert$ on the space of possible perturbations $\{\VV\}$
(the tangent space of $\VQ$ at $\Vq$),
and to restrict the search of a steepest descent direction to directions
$\VV$ of length $\Vert \VV \Vert = 1$. A natural choice would be to select
the $W^{1,\infty}$-norm\footnote{The norm $\Vert \VV\Vert_{1,\infty}$ is
the maximum between the essential supremum of $\VV$ and of its derivative
$\VD\VV$.}, with respect to which a minimizer has been shown to exist for most functionals
\cite[Prop. 3.1]{PaWeFa18}. However, in practice it is more convenient to employ
a norm that is induced by an inner product, so that the steepest descent direction
corresponds to the Riesz representative of $\rm{dJ}$ with respect to the inner product
\cite[p. 98]{HiPiUlUl09}.
\end{remark}

\subsection{Geometric transformations, moving meshes,
and parametric finite elements}
\label{sec:isoFEM}

To solve a PDE-constrained optimization problem iteratively,
it is necessary to employ a numerical method that is capable
of solving the constraint $\VA(\Vq,\Vu) = \mathbf{0}$ for any
feasible control $\Vq$. In shape optimization, this translates into
the requirement of a numerical scheme that can solve a PDE on
a domain that changes at each iteration of the optimization algorithm.
There are several competing approaches to construct a numerical
scheme with this feature, such as the level-set method \cite{AlJoTo02}
or the phase field approach \cite{BlGaFaHaSt14}, among others.

\begin{figure}[htb!]
\includegraphics[width=0.27\linewidth]{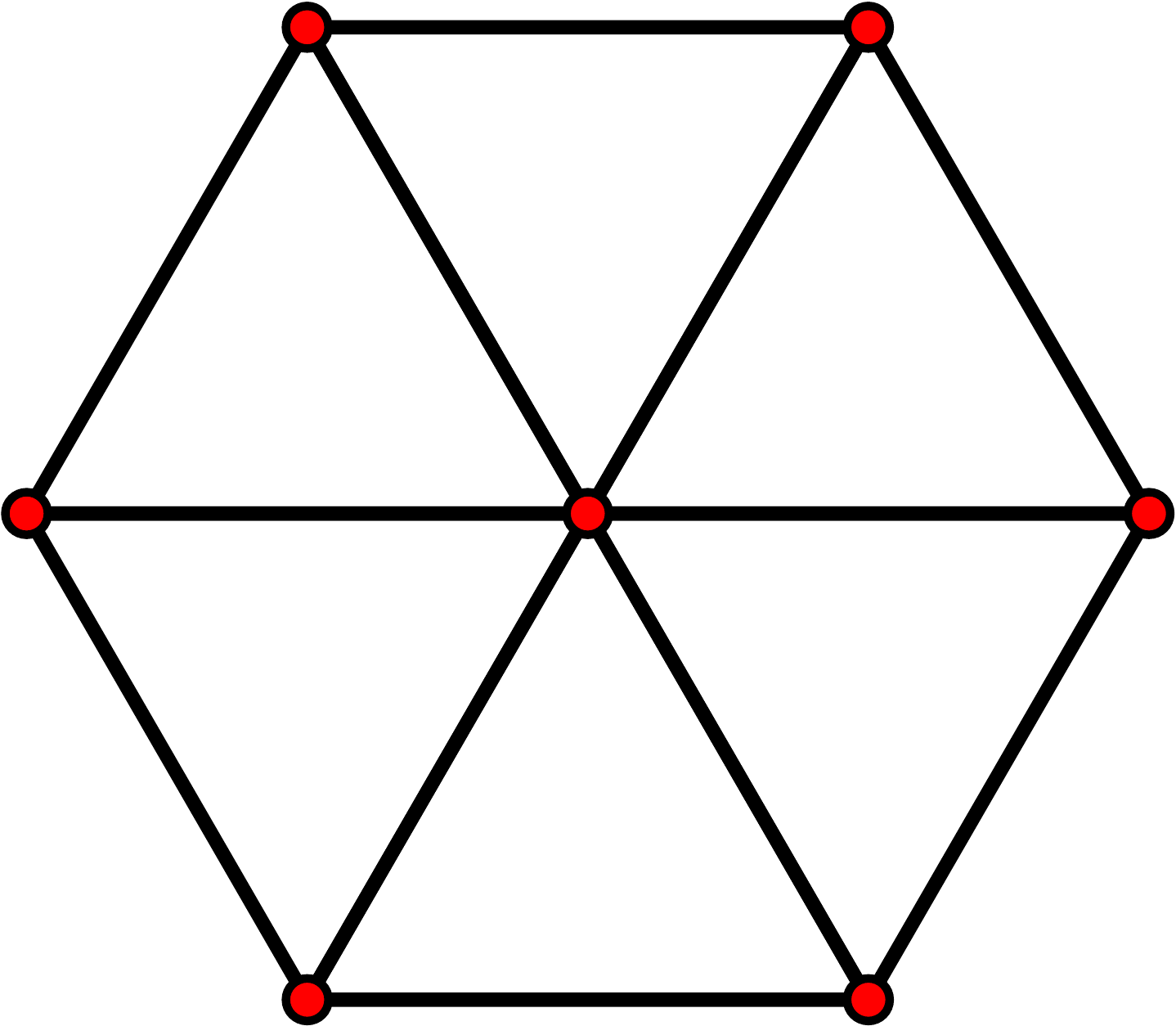}
\hspace{0.2cm}
\includegraphics[width=0.3\linewidth]{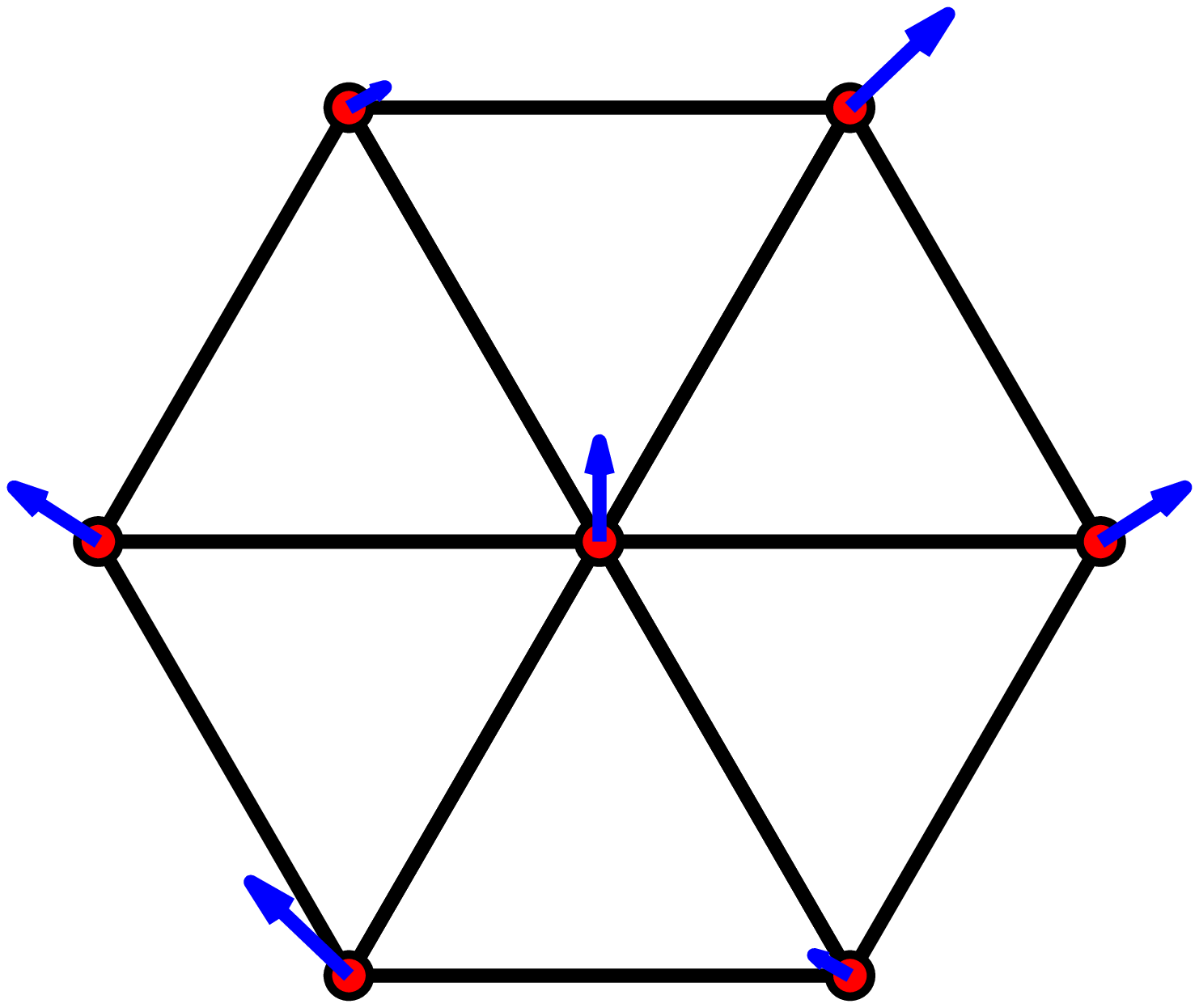}
\hfill
\includegraphics[width=0.31\linewidth]{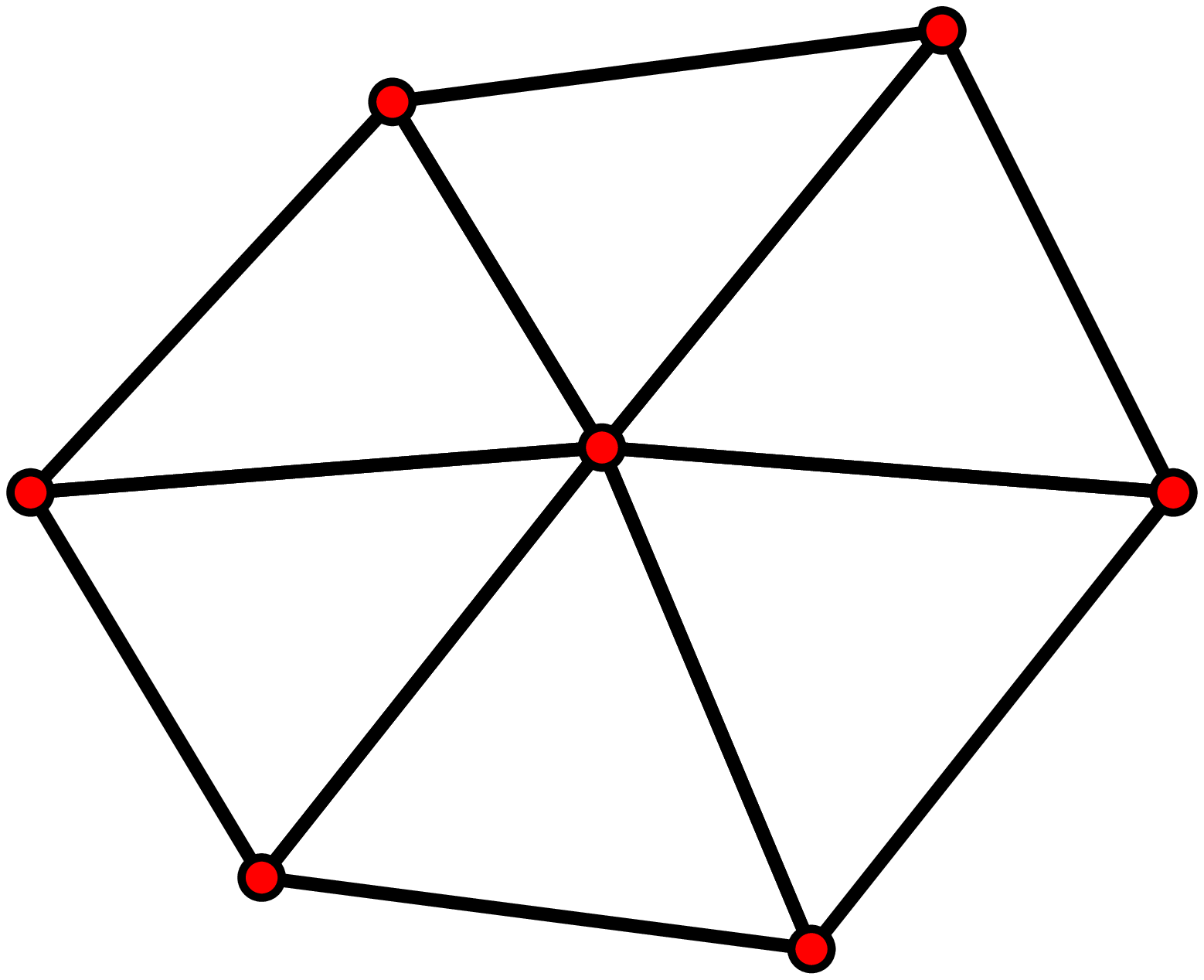}
\caption{Simplest approach to PDE-constrained shape optimization.
An initial mesh (\emph{left}) is used to compute state and adjoint variables.
This information, together with shape derivatives, is used to devise
an update of the nodes' coordinates (\emph{center}). A new and improved
initial guess is obtained by updating the nodes' coordinates and retaining
the initial mesh connectivity (\emph{right}). This process is repeated until
convergence.}
\label{fig:naivemethod}
\end{figure}

In Fireshape, we employ the approach sometimes know as 
``moving mesh'' method \cite{PaWeFa18, Al07}. In its simplest
version (see Figure \ref{fig:naivemethod}), this method replaces
(or approximates) the initial domain
$\rm{\Omega}$ with a polygonal mesh $\rm{\Omega_h}$. On this mesh,
the state and adjoint equations are solved with a standard
finite element method, whose construction on polygonal meshes is
immediate. For instance, depending on the nature of
the state constraint, one may solve the state and adjoint equations
using linear Lagrangian or P2-P1 Taylor-Hood finite elements
(see Section \ref{sec:step2}). With the state and adjoint solutions at hand,
one employs shape derivatives
(see Remark \ref{rmk:metric} and Section \ref{sec:innerproduct})
to update the coordinates of mesh nodes
while retaining the mesh connectivity of $\rm \Omega_h$.
This leads to a new mesh that represents an improved design.
This update process is repeated until some prescribed convergence
criteria are met.

The moving mesh method we just described is a simple and yet powerful method.
However, in its current formulation, it requires polyhedral meshes, which limits
the search space $\VQ$ to polyhedra. In the remaining part of this section,
we describe an equivalent interpretation of the moving mesh method
that generalizes to curved domains. Additionally, this alternative interpretation
allows approximating state and adjoint variables with arbitrarily high-order finite
elements without suffering from reduction of
convergence order due to poor approximation of domain boundaries \cite[Ch. 4.4]{Ci02}.

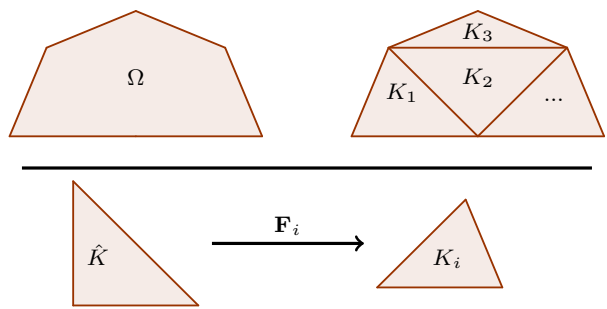
\begin{figure}[htb!]
\centering
\begin{tikzpicture}
\node at (0,0.25) {\includegraphics[width=0.4\linewidth]{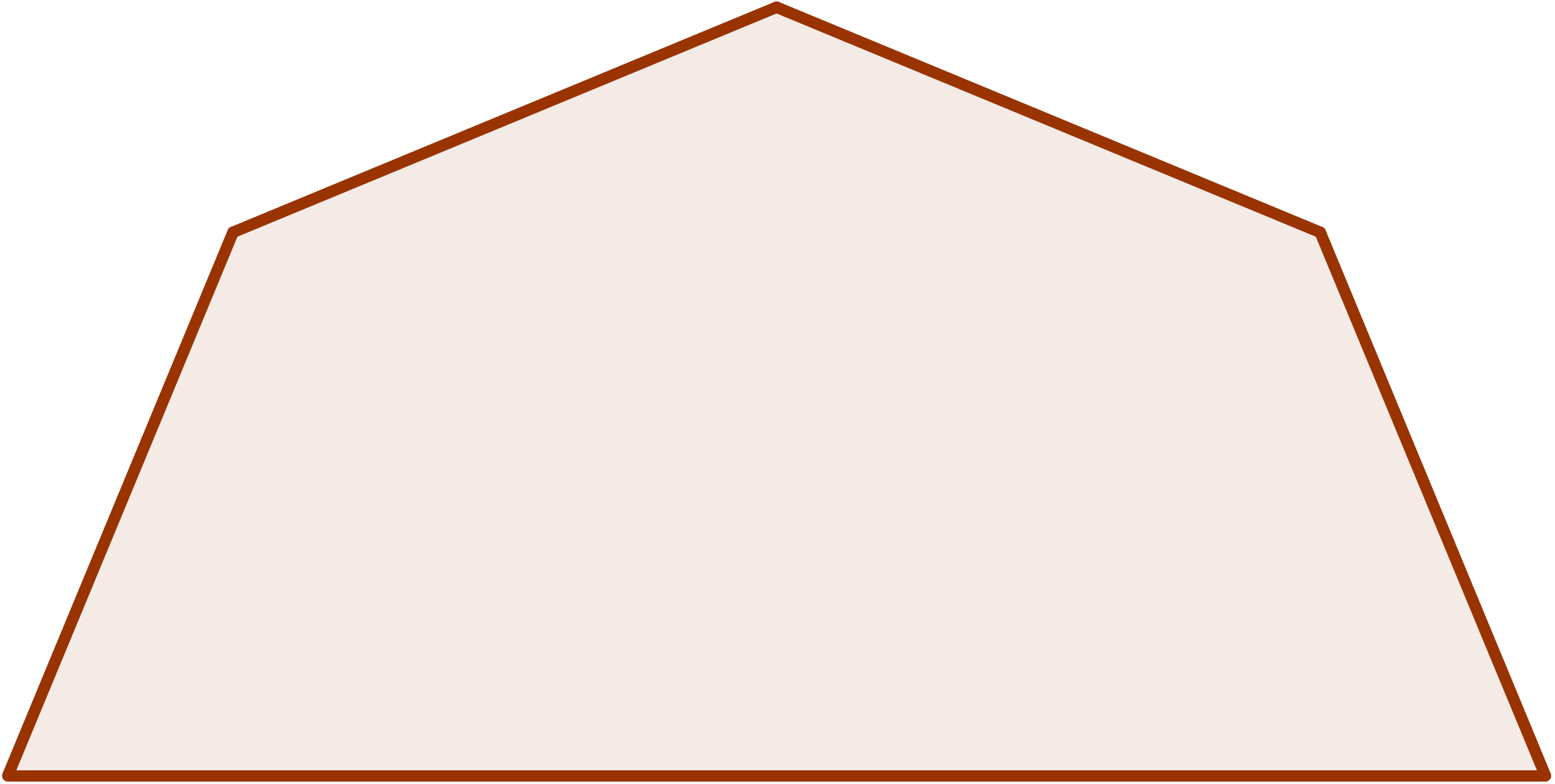}};
\node at (0, 0.2) {$\rm{\Omega}$};
\draw[very thick] (-1.5, -1) -- (6, -1);
\node at (4.5, 0.25) {\includegraphics[width=0.4\linewidth]{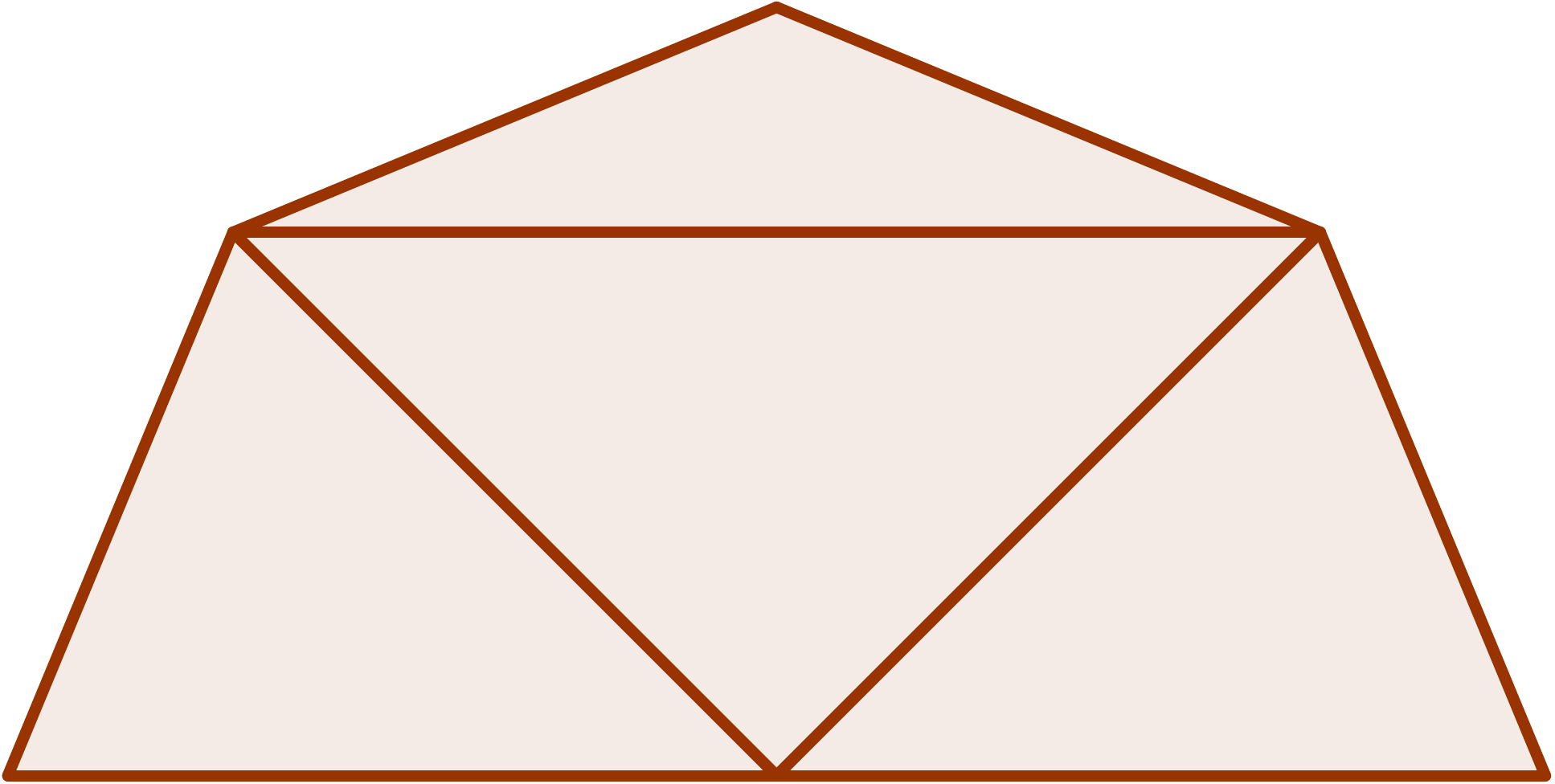}};
\node at (4.5, 0.2) {$K_2$};
\node at (3.5, 0) {$K_1$};
\node at (5.5, -0.1) {$...$};
\node at (4.5, 0.8) {$K_3$};
\node at (0,-2) {\includegraphics[width=0.2\linewidth]{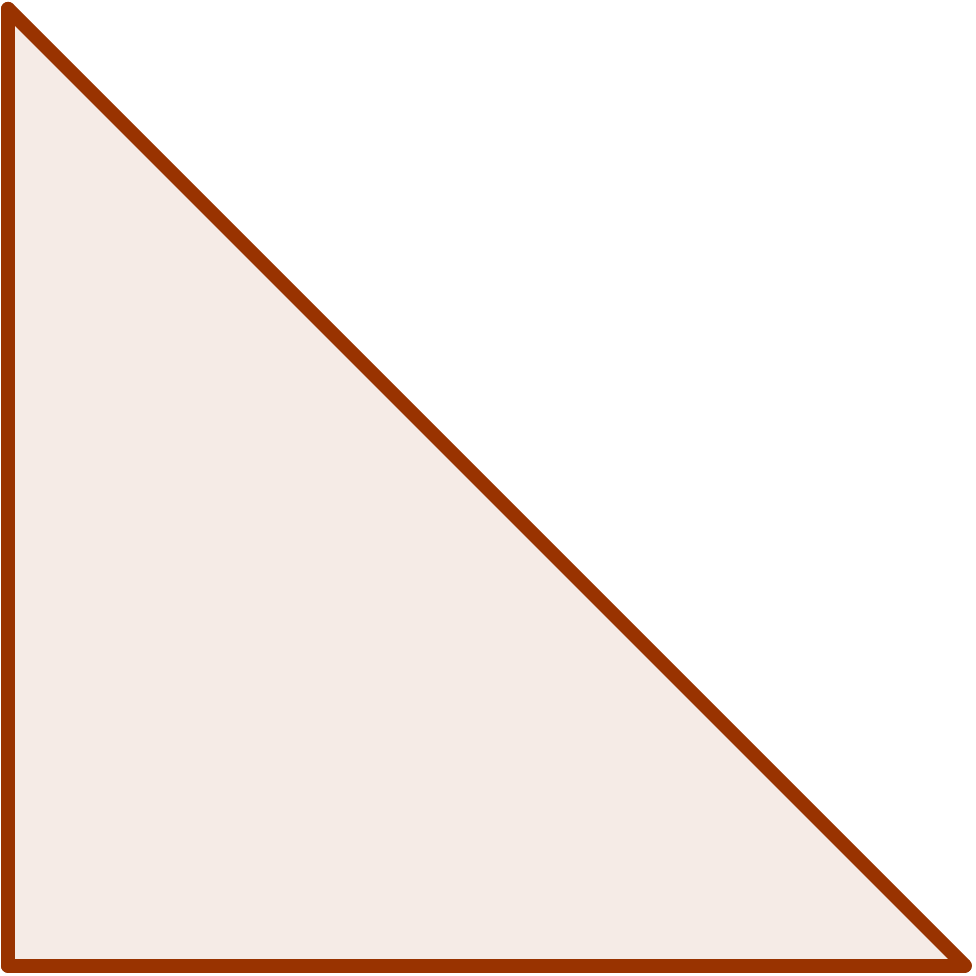}};
\node at (4,-2) {\includegraphics[width=0.2\linewidth]{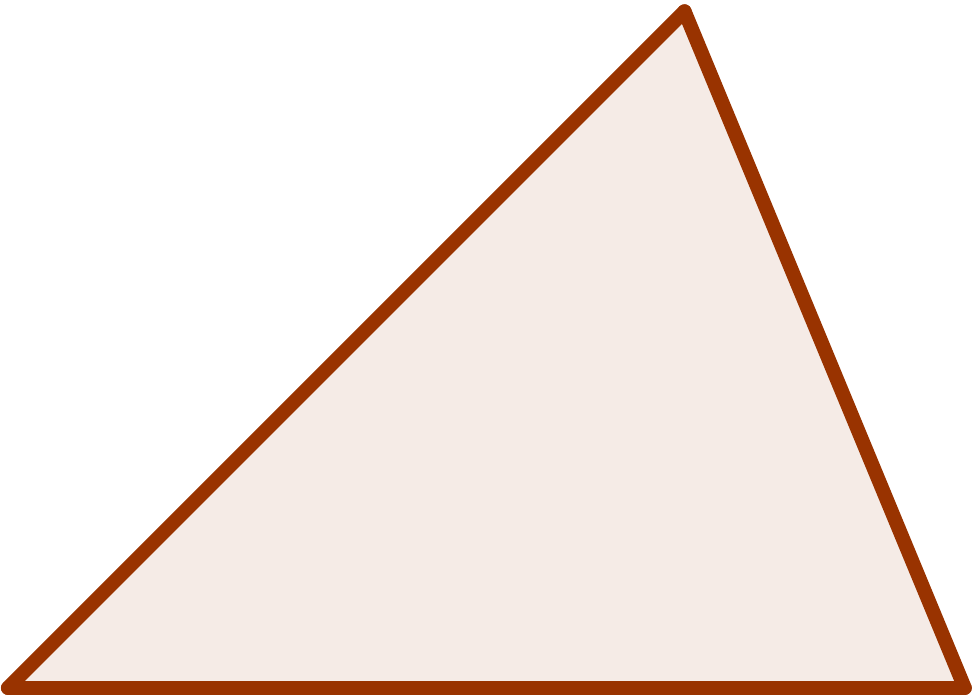}};
\draw[->, very thick] (1,-2)--(3,-2);
\node[anchor=south] at (2,-2) {$\VF_i$};
\node at (-0.5,-2.15) {$\hat{K}$};
\node at (4.1,-2.2) {$K_i$};
\end{tikzpicture}
\caption{To construct finite element basis functions, one usually begins by
triangulating a domain $\rm{Omega}$ (\emph{first row})
and then introducing diffeomorphisms
$\VF_i$ that map the reference element $\hat{K}$ to the various $K_i$
in the triangulation (\emph{second row}).}
\label{fig:FEconstruction}
\end{figure}

We begin by recalling the standard construction of parametric finite elements.
For more details, we refer to \cite[Sect. 2.3 and 4.3]{Ci02}.
To construct a parametric finite element space,
one begins by partitioning the domain $\rm{\Omega}$
into simpler geometric elements
$\{K_i\}$ (usually triangles or tetrahedra, as depicted in
the first row of Figure \ref{fig:FEconstruction}).
Then, one introduces a reference element $\hat{K}$
and a collection of diffeomorphisms $\VF_i$
that map the reference element $\hat{K}$ to the various $K_i$s, that is,
$\VF_i(\hat{K}) = K_i$ for every value of the index $i$
(as depicted in the second row of Figure \ref{fig:FEconstruction}).
Finally, one considers a set of local reference basis functions $\{\hat{b}_j\}$
defined on the reference element $\hat{K}$, and defines local
basis functions $\{b_j^i\}$ on each $K_i$ via the pullback
$b_j^i(\Vx)\coloneqq \hat{b}_j(\VF_i^{-1}(\Vx))$ for every $\Vx$ in $K_i$.
These local basis functions are used to construct global basis functions
that span the finite element space.
An important property of this construction is that the diffeomorphisms
$\VF_i$ can be expressed in terms of local linear Lagrangian
basis functions\footnote{For linear Lagrangian finite
elements, the two set of reference local basis functions $\{\hat{b}_j\}$
and $\{\hat{\beta}_m\}$ coincide.}
 $\{\hat{\beta}_m\}$, that is,
there are some coefficient vectors $\{\mubf_m^i\}$ such that
\begin{equation}\label{eq:FEdiffeo}
\VF_i(\Vx) = \sum_m \mubf_m^i \hat{\beta}_m(\Vx)
\quad \text{for every} \Vx \text{ in } \hat{K}\,.
\end{equation}
Finite elements constructed following this procedure are usually called
parametric, because they rely on the parametrization $\{\VF_i\}$.
Note that the most common finite elements families, such as Lagrangian, Raviart--Thomas, or Nedelec finite elements, are indeed parametric.

Keeping this knowledge about parametric finite elements in mind,
we can revisit the mesh method (see Figure \ref{fig:naivemethod}).
There, the main idea was to update only
nodes' coordinates and keep the mesh connectivity unchanged,
so that constructing finite elements on the new mesh is straightforward.
In \cite{PaWeFa18}, it has been shown that
the new finite element space can also be obtained
by modifying the parametric construction of finite elements on the original domain
$\rm{\Omega}$.
In the next paragraphs, we give an extended explanation (with adapted notation)
of the demonstration given in \cite{PaWeFa18}.

\begin{figure}[htb!]
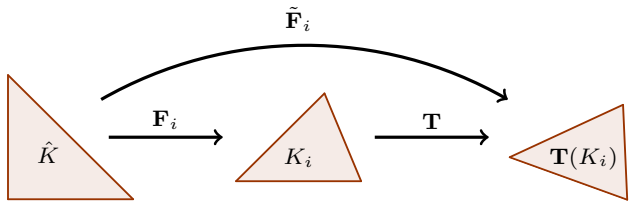

\centering
\begin{tikzpicture}
\node at (-4,-3) {\includegraphics[width=0.2\linewidth]{images/croppedtrianglereference}};
\node at (-1,-3) {\includegraphics[width=0.2\linewidth]{images/croppedtriangle}};
\node[rotate=-20] at (2.75,-3) {\includegraphics[width=0.2\linewidth]{images/croppedtriangle}};
\draw[->, very thick] (-3.5,-3)--(-2,-3);
\draw[->, very thick] (0,-3)--(1.5,-3);
\draw[->, very thick, overlay] (-3.6,-2.5) arc (120:60:5.35);
\node[anchor=south] at (-1,-1.7) {$\tilde\VF_i$};
\node[anchor=south] at (-2.75,-3) {$\VF_i$};
\node[anchor=south] at (0.75,-3) {$\VT$};
\node at (-4.3,-3.2) {{$\hat{K}$}};
\node at (-1,-3.275) {{$K_i$}};
\node at (2.75,-3.275) {{$\VT(K_i)$}};
\end{tikzpicture}
\caption{The standard construction of finite elements on the perturbed mesh
$\{\VT(K_i)\}$ employs diffeomorphisms $\tilde\VF_i$ that map the reference
element $\hat{K}$ to the various perturbed elements $\VT(K_i)$.}
\label{fig:FEconstruction_moved}
\end{figure}

Let $\VT$ denote the transformation 
employed to modify
the mesh $\rm\Omega_h$ on the left in Figure \ref{fig:naivemethod} into the
new and perturbed mesh $\VT(\rm\Omega_h)$ on the right in Figure \ref{fig:naivemethod}.
Additionally, let $\{\VT(\VK_i)\}$ denote the simple geometric
elements that constitute the latter.
Using the parametric approach,
we can construct finite elements on the new mesh
by introducing a collection of diffeomorphisms $\tilde\VF_i$
that map the reference element $\hat{K}$ to the various $\VT(K_i)$s, that is,
$\tilde\VF_i(\hat{K}) = \VT(K_i)$ for every value of the index $i$;
see Figure \ref{fig:FEconstruction_moved}.

The behavior of the transformation $\VT$ in Figure \ref{fig:naivemethod}
is prescribed only the mesh nodes.
Since its behavior on the interior of the mesh triangles can be chosen
arbitrarily, we can decide that $\VT$ is piecewise affine on each triangle.
This convenient choice implies that $\VT$ can be written
as a linear combination of piecewise affine Lagrangian finite elements
defined on the first mesh, that is,
$\VT(\Vx) = \sum_\ell \nubf_\ell B_\ell(\Vx)$ for every $\Vx$ in $\rm{\Omega}$,
where $\{\nubf_\ell\}$ are some coefficient vectors and
$\{B_\ell\}$ are global basis functions of the space of
piecewise affine Lagrangian finite elements defined on the partition $\{K_i\}$.
Since Lagrangian finite elements are constructed via pullbacks to the reference element,
for every element $K_i$ there are coefficient vectors $\{\nubf_m^i\}$ so that
the restriction $\VT\vert_{K_i}$ of $\VT$ on $K_i$ can be rewritten as
\begin{equation*}
\VT\vert_{K_i}(\Vx) = \sum_\ell \nubf_\ell B_\ell\vert_{K_i}(\Vx)
= \sum_m \tilde\nubf_m^i \hat{\beta}_m(\VF_i^{-1}(\Vx))\,.
\end{equation*}
Therefore, the composition $\VT\vert_{K_i}\circ \VF_i$ is of the form
\begin{equation}\label{eq:FEdiffeoNew}
\VT\vert_{K_i}\circ\VF_i(\Vx) = \sum_m \tilde\nubf_m^i \hat{\beta}_m(\Vx)
\quad \text{for every }\Vx\in\hat{K}\,,
\end{equation}
that is, of the same form of Equation \eqref{eq:FEdiffeo}.
This implies that, to construct finite elements on the perturbed geometry
$\VT(\rm{\Omega}_h)$, we only need to replace the original coefficients $\{\mubf_m^i\}$
in Equation \eqref{eq:FEdiffeo} with the new coefficients $\{\tilde\nubf_m^i\}$
from Equation \eqref{eq:FEdiffeoNew}.

This alternative and equivalent viewpoint on the moving mesh method
generalizes naturally to higher-order finite element approximations
{\color{black} of control and state (and adjoint) variables}. Indeed, one of
the key steps to ensure that higher-order finite elements achieve
higher-order convergence on curved domains is to employ sufficiently accurate
polynomial interpolation of domain boundaries \cite[Ch. 4.4]{Ci02}. This
boundary interpolation can be encoded in the diffeomorphisms $\{\VF_i\}$ by
using higher-order Lagrangian local basis functions. Therefore, simply
employing higher-order Lagrangian finite element transformations $\VT$ leads
to a natural extension of moving mesh method to higher-order finite elements.

This alternative and equivalent viewpoint generalizes further
to allow the use of any arbitrary discretization
of the transformation $\VT$ (for instance, using B-splines \cite{Ho03}, harmonic polynomials,
or radial basis functions \cite{We05}). The only requirement to be fulfilled to ensure the desired
order of convergence $p$ is that the maps $\VT\circ\VF_i:\Vx \mapsto \VT(\VF_i(\Vx))$ satisfy
the asymptotic algebraic estimates
\begin{equation*}
\Vert \VD^\alpha (\VT\circ\VF_i) \Vert = \Co(h^{\alpha})
\quad \text{for } 0\leq\alpha\leq p\,,
\end{equation*}
where $\VD^\alpha(\VT\circ\VF_i)$ denotes the $\alpha^{\rm{th}}$ derivative of $\VT\circ\VF_i$. Using a different discretization of the transformation $\VT$
can give several advantages, like increasing the smoothness $\VT$ (because finite elements
are generally only Lipschitz continuous) or varying how shape updates are computed
during the optimization process \cite{EiSt18}.

\begin{remark} \label{rmk:detDTsmall}

An issue that can arise with the moving mesh method is that it can lead to
poor quality (or even tangled) meshes. In terms of geometric transformations,
a mesh with poor quality corresponds to a transformation $\VT$ for which the
value $\max_{\alpha}\Vert \VD^\alpha\VT \Vert$ is large (and a tangled mesh
to a transformation $\VT$ that is not a diffeomorphism). To a certain extent,
it is possible to enforce moderate derivatives by employing suitable metrics
to extract descent directions from shape derivatives (for instance, by using
linear elasticity based inner products with a Cauchy-Riemann augmentation
\cite{IgStWe18}) and/or by adding penalty terms to the functional $\Jmat$.
\color{black}{For example, in Sections \ref{sec:finalstep} and
\ref{sec:finalstepCantilever} we add a penalization term based on the
spectral norm $\Vert\VD(\VT-\VI)\Vert_2$
of the Jacobian matrix $\VD(\VT-\VI)$. We point out that
if $\rm\Omega$ is 
$\zeta$-quasiconvex\footnote{A domain $\rm\Omega$ is $\zeta$-quasiconvex if,
for any $\Vx\in\rm\Omega$ and $\Vy\in\rm\Omega$, the length $\ell$ of the
shortest path in $\rm\Omega$ that connects $\Vx$ to $\Vy$ satisfies $\ell\leq
\zeta \Vert \Vx - \Vy\Vert$.} and $\Vert\VD(\VT-\VI)\Vert_2< \zeta$, then
$\VT:{\rm\Omega}\to\bbR^d$ is bi-Lipschitz, and thus injective.}
\end{remark}

\section{Anatomy of Fireshape}
\label{sec:anatomy}
In this section, we give more details about Fireshape's implementation and features.
Fireshape is organized in
a few core Python classes (and associated subclasses) that
implement the control space $\VQ$, the metric to be employed
by the optimization algorithm, and the (possibly PDE-constrained) objective function $\rm{J}$.
The following subsections describe these classes.

\subsection{The class \texttt{ControlSpace}}
\label{sec:controlspace}

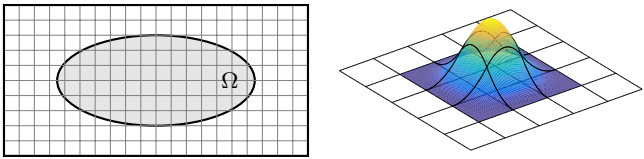
\begin{figure}[htb!]
\begin{tikzpicture}
\draw[yshift = -0.8cm,thick,fill=gray!20] (-3,1) ellipse (1.3cm and 0.6 cm);
\draw[yshift = -0.8cm,help lines, step = 0.2cm] (-5,0) grid (-1,2);
\draw[yshift = -0.8cm,thick] (-5,0) rectangle (-1,2);
\node[yshift = -0.8cm,anchor = east] at (-1.8,1) {$\rm{\Omega}$};
\node[yshift = -0.8cm,] at (1.4,1) {\includegraphics[width=4cm]{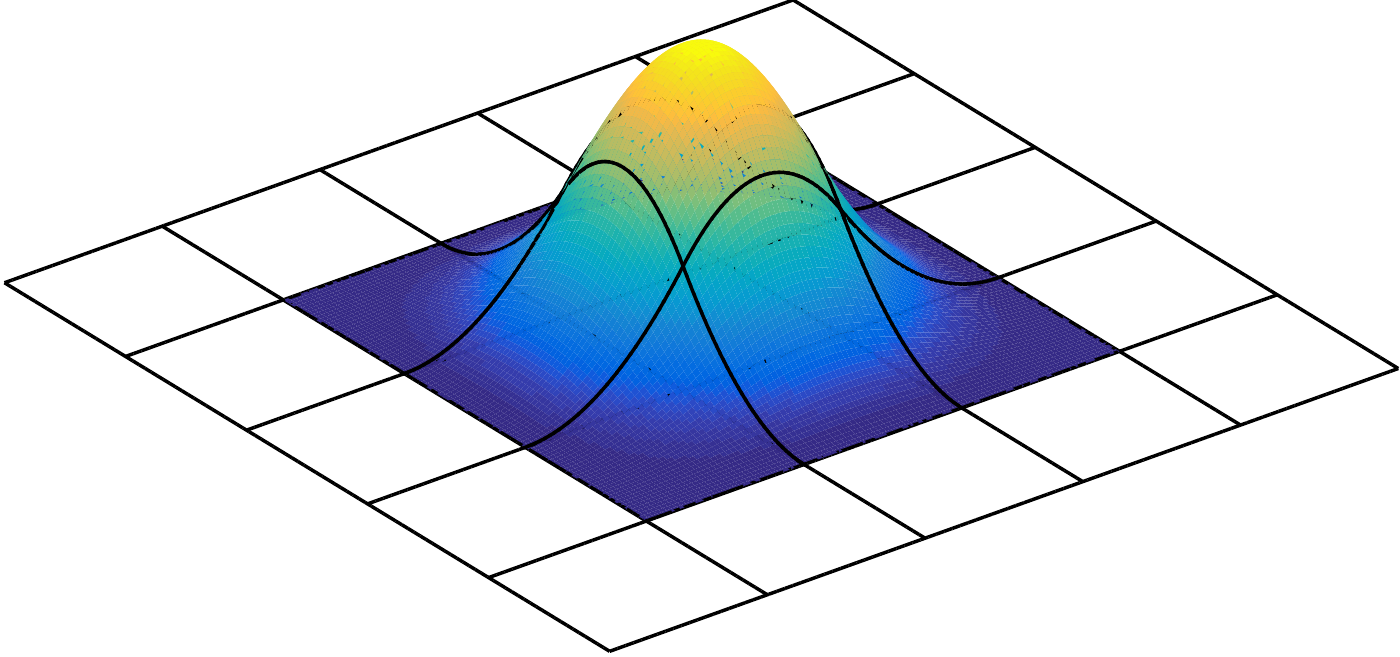}};
\end{tikzpicture}
\caption{Fireshape allows discretizing transformations $\VT$ using
tensorized B-splines. On the left, we display a Cartesian grid that
that covers the computational domain $\rm{\Omega}$. On the right,
we display a quadratic tensorized B-spline.}
\label{fig:Bsplines}
\end{figure}

Fireshape models the control space $\VQ$ of admissible domains using
geometric transformations $\VT$ as in Equation \eqref{eq:Q} (see also Figure
\ref{fig:Q}). From a theoretical perspective, the transformations $\VT$ can
be discretized in numerous different ways as long as certain minimal
regularity requirements are met. In Fireshape, the class
\texttt{ControlSpace} allows the following options: (i) Lagrangian finite
elements defined on the same mesh employed to solve the state equation, (ii)
Lagrangian finite elements defined on a mesh coarser than the mesh employed
to solve the state equation, and (iii) tensorized B-splines defined on a
separate Cartesian grid (not to be confused with a spline or B\'{e}zier
parametrization of the domain boundary, see Figure \ref{fig:Bsplines}).
{\color{black} In this case, the user must specify the boundary and the
refinement level of the Cartesian grid (\texttt{bbox} and \texttt{level}), as
well as the \texttt{order} of the underlying univariate B-splines. The user
can also set the regulariy of B-splines on the boundary of the Cartesian grid
(\texttt{voundary\_regularities}), or restrict pertubations to specific
directions (\texttt{fixed\_dims}).}

If discretization (i) can be considered to be the default option,
discretization (ii) allows introducing a regularization by discretizing
the geometry more coarsely (so-called ``regularization by discretization''),
whereas B-splines allow constructing transformation with higher regularity
(Lagrangian finite elements are only Lipschitz continuous,
whereas B-splines can be continuously differentiable and more).

The class \texttt{ControlSpace} can be easily extended to include
additional discretization options, such as radial basis functions
and harmonic polynomials.

\begin{remark}

\color{black} The current implementation of the control space $\VQ$ includes
geometric transformations that can leave the value of a shape function
unchanged. For instance, if $\Jmat(\rm{\Omega}) = \int_{\rm{\Omega}} \DX$,
then any diffeomorphism $\VT$ of the form $\VT = \VI + \VV$ with $\VV\neq
\mathbf{0}$ in $\rm{\Omega}$ and $\VV=\mathbf{0}$ on $\rm{\partial \Omega}$
satisfy $\VT(\rm{\Omega}) = \rm{\Omega}$, and thus $\Jmat(\VT(\rm{\Omega})) =
\Jmat(\rm{\Omega})$.

The presence of a shape function kernel in the control space $\VQ$ does not
represent a problem for steepest descent and quasi Newton optimization
algorithms. However, this kernel leads to singular Hessians, and Newton's
method cannot be applied in a straightforward fashion. The development of
control space descriptions that exclude shape function kernels and allow
direct applications of shape Newton methods in the framework of finite
element simulations is a highly active research field \cite{sturm2016convergence, paganini2019weakly, EtHeLoWa20}. We
predict these new control spaces will be included in Fireshape as soon as
their description has been established.

\end{remark}

\subsection{The class \texttt{InnerProduct}}
\label{sec:innerproduct}
To formulate a continuous optimization algorithm,
we need to specify how to compute lengths of vectors
(and, possibly, how to compute the angle between two vectors).
The class \texttt{InnerProduct} addresses this requirement
and allows selecting an inner product $(\cdot, \cdot)_\Ch$ to endow the control
space $\VQ$ with\footnote{Note that specifying a norm would
suffice to formulate an optimization algorithm. However,
as mention in Remark \ref{rmk:metric}, it is computationally
more convenient to restrict computations to inner product spaces.}.

The choice of the inner product affects how steepest-descent directions are
computed. Indeed, a steepest-descent direction is a direction $\VV$ of length
$\Vert \VV\Vert_\Ch = 1$ such that $\rm{dJ}(\Vq,\VV)$ is minimal. Let
$\alpha\coloneqq \Vert \rm{dJ}(\Vq,\cdot)\Vert_*$ denote the length of the
operator {\color{black}$ \rm{dJ}(\Vq,\cdot)$} measured with respect to the dual
norm. Then \cite[p. 103]{HiPiUlUl09}, the steepest descent direction $\VV$
satisfies the equation \begin{equation*} \alpha(\VV,\VW)_\Ch =
-\rm{dJ}(\Vq,\VW)\quad \text{for all } \VW \text{ in } \Ch\,, \end{equation*}
which clearly depends on the inner product $(\cdot, \cdot)_\Ch$.

The control space $\VQ$ can be endowed with different inner products.
In Fireshape, the class \texttt{InnerProduct} allows the following
options: (i) an $H^1(\rm{\Omega})$ inner product based on standard
Galerkin stiffness and mass matrices, (ii) a Laplace inner product
based on the Galerkin stiffness matrix, and (iii) an elasticity inner product
based on the linear elasticity mechanical model. These options can
be complemented with additional homogeneous Dirichlet boundary
conditions to specify parts of the boundary $\partial\rm{\Omega}$
that are not to be modified during the shape optimization procedure.

Although all three options are equivalent from a theoretical perspective,
in practice it has been observed that option (iii) generally leads to
geometry updates that result in meshes of higher equality compared to
options (i) and (ii). A thorough comparison is available in \cite{IgStWe18},
where the authors also suggest to consider complementing these inner products
with terms stemming from Cauchy-Riemann equations to further increase
mesh quality. This additional option is readily available in Fireshape.

\subsection{The classes \texttt{Objective} and \texttt{PdeConstraint}}
\label{sec:objective}

In the vast majority of cases, users who aim to solve a PDE-constrained
shape optimization problem are only required to instantiate
the two classes \texttt{Objective} and \texttt{PdeConstraint}, where
they can specify the formula of the function $\rm{J}$ to be minimized
and the weak formulation of its PDE-constraint
$\VA(\rm{\Omega}, \Vu) = \bold{0}$ (see Sections \ref{sec:step2} and \ref{sec:step3}, for instance).
Since Fireshape is built on top
of the finite element library Firedrake, these formulas must be
written using the Unified Form Language (UFL). We refer to
the tutorials on the website \cite{Firedrakewebsite}
for more details about Firedrake and UFL.

\subsection{Supplementary classes}
\label{sec:extraclasses}

Fireshape also includes a few extra classes to specify
additional constraints, such as volume or perimeter
constraints on the domain $\rm{\Omega}$, or
spectral constraints to control the singular values of
the transformation $\VT$. For more details about these extra options,
we refer to Fireshape's documentation and tutorials \cite{Fireshapedocu}.

\section{Conclusions}
\label{sec:conclusions}
We have introduced Fireshape: an open-source and automated
shape optimization toolbox for the
finite element software Firedrake. Fireshape is based on the moving mesh
method and allows users with minimal shape optimization knowledge
to tackle challenging shape optimization problems constrained to 
PDEs. In particular, Fireshape computes adjoint equations
and shape derivatives in an automated fashion, allows decoupled discretizations
of control and state variables, and gives full access to Firedrake's and PETSc's
discretization and solver capabilities as well as to the extensive optimization library
ROL.

\begin{acknowledgements}
The work of Florian Wechsung was partly supported by the EPSRC
Centre For Doctoral Training in Industrially Focused Mathematical Modelling
[grant number EP/L015803/1].\\~\\
\textbf{Contributions}
Alberto Paganini and Florian Wechsung have contributed equally
to the development of the software Fireshape and to this manuscript.\\~\\
\textbf{Conflict of interest statement}
 On behalf of all authors, the corresponding author states that there is no conflict of interest.
 \\~\\
\textbf{Replication of results}
For reproducibility, we cite archives of the exact software versions
used to produce the results in this paper. All major Firedrake
components have been archived on Zenodo~\cite{zenodo-firedrake-and-driver}. An
installation of Firedrake with components matching those used to
produce the results in this paper can by obtained following the
instructions at
\begin{center}
\url{https://www.firedrakeproject.org/download.html} .
\end{center}
The exact version of the Fireshape library used for these results has also been archived at~\cite{zenodo-fireshape}.
The latest version of the Fireshape library can be found at 
\begin{center}
    \url{https://github.com/Fireshape/Fireshape} .
\end{center}
\end{acknowledgements}

\bibliographystyle{plainurl}

\end{document}

%% file: slidesymbols.tex
\providecommand{\Div}{\operatorname{div}}          

\renewcommand{\Re}{\operatorname{Re}}             

\newcommand{\Va}{{\mathbf{a}}}

\newcommand{\Vg}{{\mathbf{g}}}

\newcommand{\Vn}{{\mathbf{n}}}

\newcommand{\Vp}{{\mathbf{p}}}
\newcommand{\Vq}{{\mathbf{q}}}

\newcommand{\Vu}{{\mathbf{u}}}
\newcommand{\Vv}{{\mathbf{v}}}

\newcommand{\Vx}{{\mathbf{x}}}
\newcommand{\Vy}{{\mathbf{y}}}

\providecommand{\Bx}{{\boldsymbol{x}}}

\newcommand{\VA}{{\mathbf{A}}}

\newcommand{\VD}{{\mathbf{D}}}

\newcommand{\VF}{{\mathbf{F}}}

\newcommand{\VI}{{\mathbf{I}}}

\newcommand{\VK}{{\mathbf{K}}}

\newcommand{\VQ}{{\mathbf{Q}}}

\newcommand{\VT}{{\mathbf{T}}}
\newcommand{\VU}{{\mathbf{U}}}
\newcommand{\VV}{{\mathbf{V}}}
\newcommand{\VW}{{\mathbf{W}}}

\newcommand{\mubf}{\boldsymbol{\mu}}
\newcommand{\nubf}{\boldsymbol{\nu}}

\newcommand{\Jmat}{{\rm J}}

\providecommand{\Ch}{{\cal H}}

\providecommand{\Cl}{{\cal L}}

\providecommand{\Co}{{\cal O}}

\providecommand{\bbN}{\mathbb{N}}

\providecommand{\bbR}{\mathbb{R}}

\newcommand{\DX}{\,\mathrm{d}\Bx}
\newcommand{\dS}{\,\mathrm{d}S}

